\documentclass{article}
\usepackage{personal}
\usepackage{mathcros}
\usepackage{algebra}
\usepackage{numbertheory}
\usepackage[all,tips,dvips]{xy}

\CompileMatrices

\newcommand\Alg{\operatorname{\mathbf{Alg}}}
\newcommand\Ab{\operatorname{\mathbf{Ab}}}
\newcommand\hexp{\operatorname{hexp}}

\numberwithin{equation}{section}

\title{The Theory of Witt Vectors}
\author{Joseph Rabinoff}

\begin{document}

\maketitle

\renewcommand\numberline[1]{\makebox[2em][l]{\large\bfseries #1.}}
\renewcommand\contentsline[4]{{\bfseries #2}\dotfill\makebox[1.5em][r]{#3}\\}
\tableofcontents

\bigskip\bigskip\noindent
The theory of Witt vectors, while mostly elementary, manages to package such
algebraic power into the letter $W$ that it turns up in many areas of
mathematics.  And whereas the Witt rings enjoy a great number of
``well-known'' properties, it can be difficult to find a written statement
of these properties, much less a proof.  
Difficult, that is, if one does
not know that almost everything one could want to know about the Witt
vectors, including the vast majority%
\footnote{Basically, everything except Section~\ref{sec:motivation}.}
of the material below, is contained in Hazewinkel's book~\cite{hazewinkel}.  
As I was unaware of this reference until B.~Conrad kindly pointed it out, 
I undertook to write a reference manual for the theory of Witt vectors.
The results I've included are simply those which I've come across in
my own research; I do not by any means claim that they are comprehensive.
My contribution to the material below mainly consists of the specific
statements and their proofs;%
\footnote{Again, as I was unaware that everything is proved in
  \cite{hazewinkel}, many of the proofs are my own, and differ from those
  in \cite{hazewinkel}.  The presentation is otherwise similar.}
none of it can be considered original work. 
If the reader finds an omission or error, I would be grateful if she would 
email me at \texttt{rabinoff@post.harvard.edu}.

As the theory of Witt vectors is quite elementary, the reader need only be
familiar with the fundamental concepts and constructions from commutative
algebra and category theory to understand the vast majority of this paper.
However, experience with $p$-adic rings and other linearly topological rings
would be helpful.

In this paper, all rings are commutative with identity element $1$ (unless
stated otherwise).
For a ring $K$, we let $\Alg_K$ denote the category of $K$-algebras.  Let
$\Ab$ denote the category of abelian groups.

\section{Motivation}
\label{sec:motivation}

The theory of Witt vectors arises from the study of the properties of
certain $p$-adic rings, and indeed, provides a construction of the unramified
extensions of the $p$-adic integers, for example.  These rings are defined
as follows:

\begin{defn}
  Let $p$ be a prime number.  A ring $R$ is called a 
  \emph{strict $p$-ring} provided that $R$ is complete and Hausdorff with
  respect to the 
  $p$-adic topology, $p$ is not a zero-divisor in $R$, and the residue
  ring $K = R/p$ is perfect (i.e., the map $x\mapsto x^p$ is bijective on
  $K$). 
\end{defn}

For our purposes, the following are the important properties of strict
$p$-rings (cf. \cite[Chapter~II]{serre}):

\begin{thm} \label{thm:strictprings}
  Let $K$ be a perfect ring of characteristic $p$.
  \begin{enumerate}
  \item There is a strict $p$-ring $R$ with residue ring $K$, which is
    unique up to canonical isomorphism.
  \item There exists a unique system of representatives $\tau: K\to R$,
    called the \emph{Teichm\"uller representatives},
    such that $\tau(xy) = \tau(x)\tau(y)$ for all $x,y\in K$.
  \item Every element $x$ of $R$ can be written uniquely in the form
    $x = \sum_{n=0}^\infty \tau(x_n)\,p^n$ for $x_n\in K$.
  \item The formation of $R$ and $\tau$ is functorial in $K$, in that if
    $f:K\to K'$ is a homomorphism of perfect rings of characteristic $p$,
    and $R'$ is the strict $p$-ring with residue ring $K'$ and section
    $\tau'$, then there is 
    a unique homomorphism $F:R\to R'$ making the following squares commute:
    \[\xymatrix{
      R \ar[d] \ar[r]^F & {R'} \ar[d] & & R \ar[r]^F & {R'}  \\
      K \ar[r]^f & {K'}  & & K \ar[u]^\tau \ar[r]^f & {K'} \ar[u]^{\tau'}
    }\]
    The map $F$ is given by
    \[ F\left(\sum_{n=0}^\infty \tau(x_n)\,p^n\right)
    = \sum_{n=0}^\infty \tau'(f(x_n))\,p^n. \]
  \end{enumerate}

\end{thm}

\begin{eg}
  Let $R$ be an unramified extension of $\Z_p$, with residue field
  $K = R/p\cong\F_q$.  Then $R$ is a strict $p$-ring, and is hence the
  unique strict $p$-ring with residue field $\F_q$.  The 
  Teichm\"uller representatives are constructed as follows: we have
  $\F_q^\times\cong\Z/(q-1)\Z$, so that the nonzero elements of $\F_q$ are
  the roots of the polynomial $X^{q-1}-1$.  By Hensel's Lemma, each
  element $x$ of $\F_q^\times$ has a unique lift $\tau(x)\in R$ also
  satisfying $\tau(x)^{q-1} - 1 = 0$.  Setting $\tau(0) = 0$ completes the
  definition of the map $\tau$.  In other words, the Teichm\"uller
  representatives are exactly the $(q-1)$st roots of unity in $R$, union
  $\{0\}$. 
\end{eg}

Theorem~\ref{thm:strictprings} is an abstract fact, which may be proved
without yielding a useful construction of the strict $p$-ring $R$ with
given residue ring $K$.  But the strong unicity and functoriality of $R$
indicates that one should be able to construct it algebraically in terms
of $K$.  We can certainly reconstruct the set underlying $R$ as all sums of
the form $\sum_{n=0}^\infty \tau(x_n)\,p^n$, where we think of $\tau(x_n)$
as a parameter depending only on $x_n\in K$.  But in order to reconstruct
the ring structure on $R$, we need to understand the addition and
multiplication laws in terms of the arithmetic of $K$.  Put another way,
if $\sum \tau(x_n)\,p^n + \sum \tau(y_n)\,p^n = \sum \tau(s_n)\,p^n$, we
need to write the $s_n$ in terms of the $x_n$ and $y_n$, and similarly for
multiplication.  By unicity, the answer should not depend on $R$, and by
functoriality, the answer should not even depend on $K$, in that the same
addition and multiplication laws will have to work for every $K$.  This
suggests that the $s_n$ will be given by polynomials in the $x_n$ and
$y_n$ with $p$-integral rational coefficients; that is, by polynomials
whose coefficients are rational numbers with nonnegative $p$-adic valuation.
In fact the $s_n$ will be given by \emph{integer} polynomials in the $x_n$
and $y_n$.  

The following lemma will fundamental in proving the $p$-integrality of
polynomials:

\begin{lem} \label{lem:pcong}
  Let $A$ be a ring, and let $x,y\in A$ be such that 
  $x\equiv y\pmod{pA}$.  Then for all $i\geq 0$ we have
  $x^{p^i}\equiv y^{p^i}\pmod {p^{i+1}A}$.
\end{lem}

\begin{pf}
  We proceed by induction on $i$; the case $i = 0$ is clear.  Let 
  $i\geq 1$, and write 
  $x^{p^{i-1}} = y^{p^{i-1}} + p^iz$ for $z\in A$.  Raising both sides
  side to the $p$th power, we obtain
  \[ x^{p^i} = y^{p^i} 
  + \sum_{n=1}^{p-1} \binom{p}{n} y^{p^i(p-n)}p^{in}z^n
  + p^{ip} z^p. \]
  The lemma follows because $p$ divides all of the binomial coefficients,
  and $ip \geq i+1$.

\end{pf}

Let $R$ be a strict $p$-ring with residue ring $K$, and suppose that
$\sum \tau(x_n)\,p^n + \sum \tau(y_n)\,p^n = \sum \tau(s_n)\,p^n$.
To calculate the $s_n$, we proceed inductively.  Looking mod $p$, we have
\[ \tau(x_0) + \tau(y_0) \equiv \tau(s_0)\pmod p, \]
so since $\tau(x) = x \pmod p$, we have $x_0 + y_0 = s_0$.  The na\"ive
second step is to write
\[ \tau(x_0) + p\tau(x_1) + \tau(y_0) + p\tau(y_1)
\equiv \tau(s_0) + p\tau(s_1) 
\equiv \tau(x_0 + y_0) + p\tau(s_1) \pmod{p^2}, \]
then rewrite to find
\[ p\tau(s_1) \equiv 
\tau(x_0) + \tau(y_0) - \tau(x_0 + y_0) + p(\tau(x_1) + \tau(y_1))
\pmod{p^2}. \]
But whereas we know that 
$\tau(x_0) + \tau(y_0) - \tau(x_0+y_0)\equiv 0\pmod p$, we have no idea
what its residue mod $p^2$ is.  The trick to calculating $s_1$ is as
follows.  Since $K$ is perfect, every $x\in K$ has a unique $p$th root,
written $x^{1/p}$.  Since $x_0 + y_0 = s_0$, we must have
$x_0^{1/p} + y_0^{1/p} = s_0^{1/p}$.  By Lemma~\ref{lem:pcong} above and
the fact that $\tau$ commutes with multiplication, we can write
\[ \tau(s_0) = \tau(s_0^{1/p})^p = \tau(x_0^{1/p} + y_0^{1/p})^p \equiv 
(\tau(x_0^{1/p}) + \tau(y_0^{1/p}))^p \pmod{p^2}. \]
Therefore,
\begin{equation*}
p\tau(s_1) \equiv 
\tau(x_0^{1/p})^p + \tau(y_0^{1/p})^p 
- (\tau(x_0^{1/p}) + \tau(y_0^{1/p}))^p + p(\tau(x_1) + \tau(y_1))
\pmod{p^2}. \end{equation*}
Expanding out the above equation and dividing by $p$, we obtain
\[ \tau(s_1) \equiv \tau(x_1) + \tau(y_1) - 
\sum_{n=1}^{p-1} \frac 1p \binom pn \tau(x_0^{n/p})\tau(y_0^{(p-n)/p})
\pmod p, \] 
and therefore,
\[ s_1 = x_1 + y_1 - \sum_{n=1}^{p-1} \frac 1p \binom pn 
x_0^{n/p} y_0^{(p-n)/p}. \]

The above bit of formal manipulation has nothing to do with $K$.
Indeed, let $X_0,X_1,Y_0,Y_1$ be indeterminates, and define polynomials
$w_1(X_0) = X_0$ and $w_p(X_0,X_1) = X_0^p + pX_1$; then solve
the polynomial equations
\[\begin{split}
 S_0 = w_1(S_0) &= w_1(X_0) + w_1(Y_0) = X_0 + Y_0 \\
 S_0^p + pS_1 = w_p(S_0,S_1) &= w_p(X_0,X_1) + w_p(Y_0,Y_1) 
 = X_0^p + pX_1 + Y_0^p + pY_1
\end{split}\]
for $S_0$ and $S_1$.  This is the same bit of algebra as above, except
with $X_0$ replacing $\tau(x_0^{1/p})$, etc., so we have
\[\begin{split}
  S_0 &= X_0 + Y_0 \\
  S_1 &= X_1 + Y_1 - \sum_{n=1}^{p-1}\frac 1p\binom pn X_0^n Y_0^{p-n}.
\end{split}\]
In particular, $S_0\in\Z[X_0,Y_0]$ and $S_1\in\Z[X_0,X_1,Y_0,Y_1]$.
Substituting $\tau(x_0^{1/p})$ back in for $X_0$, etc., we see that 
$s_0 = S_0(x_0,y_0)$ and $s_1 = S_1(x_0^{1/p},y_0^{1/p},x_1,y_1)$.

Witt \cite{witt} realized this, and also discovered the pattern.  Define 
\[ w_{p^n}(X_0,X_1,\ldots,X_n) = \sum_{i=0}^n p^i X_i^{p^{n-i}}
= X_0^{p^n} + pX_1^{p^{n-1}} + \cdots + p^{n-1}X_{n-1}^p + p^nX_n. \]
Letting $X_1,Y_1,X_2,Y_2,\ldots$ be indeterminates, inductively
find $S_n$ that solve the polynomial equations
\begin{equation}\label{eq:defCi}
  w_{p^n}(S_0,S_1,\ldots,S_n) = w_{p^n}(X_0,X_1,\ldots,X_n) + 
  w_{p^n}(Y_0,Y_1,\ldots,Y_n). 
\end{equation}
As the only term of $w_{p^n}(S_0,S_1,\ldots,S_n)$ involving $S_n$ is
$p^nS_n$, it is clear that there are unique polynomials $S_n$ with
rational coefficients satisfying the above identity.  What Witt showed is
that in fact, $S_n \in \Z[X_0,Y_0,X_1,Y_1,\ldots,X_n,Y_n]$.  Assuming
this, we have:

\begin{thm} \label{thm:wittpadd}
  Let $R$ be a strict $p$-ring, let $K = R/p$ be its residue ring, and let
  $\tau:K\to R$ the system of Teichm\"uller representatives.  Suppose that
  \[ \sum_{n=0}^\infty \tau(x_n)\,p^n + \sum_{n=0}^\infty \tau(y_n)\,p^n 
  = \sum_{n=0}^\infty \tau(s_n)\,p^n. \]
  Then with the $S_n$ as above, we have
  \[ s_n = S_n\big(x_0^{1/p^n},y_0^{1/p^n},x_1^{1/p^{n-1}},y_1^{1/p^{n-1}},
  \ldots,x_{n-1}^{1/p},y_{n-1}^{1/p},x_n,y_n\big). \]
\end{thm}

\begin{pf*}
  Let $n\geq 0$, and define $\td s_i$ to be 
  \[ \td s_i := S_{p^i}(x_0^{1/p^i},y_0^{1/p^i},x_1^{1/p^{i-1}},y_1^{1/p^{i-1}},\ldots,x_i,y_i) \]
  for $i \leq n$.  
  This is a polynomial identity in $K$ with integer coefficients, so we
  may take $p^{n-i}$th roots to obtain
  \[ \td s_i^{1/p^{n-i}} = S_{p^i}(x_0^{1/p^n},y_0^{1/p^n},x_1^{1/p^{n-1}},y_1^{1/p^{n-1}},\ldots,x_i^{1/p^{n-i}},y_i^{1/p^{n-i}}). \]
  Of course this is the same as saying that
  \[\begin{split}
    &S_{p^i}\big(\tau(x_0^{1/p^n}),\tau(y_0^{1/p^n}),\tau(x_1^{1/p^{n-1}}),\tau(y_1^{1/p^{n-1}}),\ldots,\tau(x_i^{1/p^{n-i}}),\tau(y_i^{1/p^{n-i}})\big) \\
    &\qquad\equiv \tau(\td s_i^{1/p^{n-i}}) \pmod p, 
  \end{split}\]
  so by Lemma~\ref{lem:pcong},
  \begin{equation}\begin{split}\label{eq:huge}
    &S_{p^i}\big(\tau(x_0^{1/p^n}),\tau(y_0^{1/p^n}),\tau(x_1^{1/p^{n-1}}),\tau(y_1^{1/p^{n-1}}),\ldots,\tau(x_i^{1/p^{n-i}}),\tau(y_i^{1/p^{n-i}})\big)^{p^{n-i}} \\
    &\qquad\equiv \tau(\td s_i^{1/p^{n-i}})^{p^{n-i}} 
    = \tau(\td s_i) \pmod{p^{n-i+1}}.
  \end{split}\end{equation}

  Substituting $\tau(x_i^{1/p^{n-i}})$ for $X_i$ and
  $\tau(y_i^{1/p^{n-i}})$ for $Y_i$, into \eqref{eq:defCi},
  we have the identity
  \begin{equation*}\begin{split}
    \tau(x_0) + & p\tau(x_1) + \cdots + p^n\tau(x_n) +
    \tau(y_0) + p\tau(y_1) + \cdots + p^n\tau(y_n) \\ &=
    S_1\big(\tau(x_0^{1/p^n}),\tau(y_0^{1/p^n})\big)^{p^n} +
    pS_p\big(\tau(x_0^{1/p^n}),\tau(y_0^{1/p^n}),\tau(x_1^{1/p^{n-1}}),\tau(y_1^{1/p^{n-1}})\big)^{p^{n-1}} \\ &\quad+\cdots+
    p^nS_{p^n}\big(\tau(x_0^{1/p^n}),\tau(y_0^{1/p^n}),\tau(x_1^{1/p^{n-1}}),\tau(y_1^{1/p^{n-1}}),\ldots,\tau(x_n),\tau(y_n)\big)
  \end{split}\end{equation*}
  in $R$.  Using \eqref{eq:huge}, we may rewrite the right-hand side of
  the above equation as
  \[\begin{split}
    \tau(x_1) + & p\tau(x_1) + \cdots + p^n\tau(x_n) +
    \tau(y_1) + p\tau(y_1) + \cdots + p^n\tau(y_n) \\ &\equiv
    \tau(\td s_1) + p\tau(\td s_1) + \cdots + p^n\tau(\td s_n) 
    \pmod{p^{n+1}}.
  \end{split}\]
  This implies that $\td s_i = s_i$ for $i\leq n$, as desired.

\end{pf*}

Witt also showed the analogous result for multiplication: namely, if we
solve the polynomial equations
\[ w_{p^n}(Z_0,Z_1,\ldots,Z_n) =
w_{p^n}(X_0,X_1,\ldots,X_n)w_{p^n}(Y_0,Y_1,\ldots,Y_n) \]
then $Z_n\in\Z[X_0,Y_0,X_1,Y_1,\ldots,X_n,Y_n]$.  Setting 
\[ z_n = Z_n\big(x_0^{1/p^n},y_0^{1/p^n},x_1^{1/p^{n-1}},y_1^{1/p^{n-1}},
\ldots,x_{n-1}^{1/p},y_{n-1}^{1/p},x_n,y_n\big), \]
the reader can verify (following the proof of Theorem~\ref{thm:wittpadd}) that
\[ \left(\sum \tau(x_n)\,p^n\right)\cdot\left(\sum\tau(y_n)\,p^n\right)
 = \sum\tau(z_n)\,p^n. \]

Hence the ring laws on $R$ are described entirely by the polynomials $S_n$
and $Z_n$, which were obtained algebraically.  This motivates the
definition of the Witt vectors: given any ring $A$, we will construct a
ring $W_p(A)$, whose addition and multiplication laws are somehow given by
the polynomials $S_n$ and $Z_n$.  Taking $A = K$, we will show 
(Theorem~\ref{thm:wittpring}) that $W_p(K)$
is the strict $p$-ring with residue ring $K$.

\section{Definition of the Witt Rings}
\label{sec:construction} 

There are several flavors of Witt rings, so for the sake of uniformity in
the statements of results, we will define Witt rings associated to the
following subsets of the natural numbers:

\begin{defn}
  A subset $P\subset\N = \{1,2,3,\ldots\}$ is a \emph{divisor-stable set}
  provided that 
  $P\neq\emptyset$, and if $n\in P$, then all proper divisors of $n$ are
  also in $P$.  If $P$ is a divisor-stable set, we let
  $\wp(P)$ denote the set of prime numbers contained in $P$.
\end{defn}

\begin{rem}
  Let $P$ be a divisor-stable set.  We make the following observations:
  \begin{enum}
  \item Since $P\neq\emptyset$, we automatically have $1\in P$.
  \item If $n\in P$, then all prime factors of $n$ are contained in
    $\wp(P)$. 
  \item The multiplicatively closed subset $T$ generated by $P$ is simply
    the set of products of primes in $\wp(P)$.
  \end{enum}

\end{rem}

\begin{eg}
  \begin{enum}
  \item The set $\N$ is divisor-stable, as are the finite sets
    $\{1, 2, \ldots, n\}$.
  
  \item Let $p$ be a prime number.  The set 
    $P_p = \{1, p, p^2, p^3, \ldots\}$ is divisor-stable, as are the
    finite sets $P_{p(n)} = \{1, p, p^2, \ldots, p^n\}$.
  \end{enum}

\end{eg}

The following is key to defining the Witt rings.

\begin{defn}
  Let $n\in\N$.  Define the \emph{$n$th Witt polynomial} to be
  \[ w_n = \sum_{d\mid n} dX_d^{n/d} \in \Z[\{X_d:d\mid n\}]. \]
  For any divisor-stable set $P$ and any ring $A$, define the set 
  \[ W_P(A) := \prod_{n\in P} A = A^P; \]
  and for $x\in W_P(A)$, we write $x_n$ for the $n$th coordinate, so
  $x = (x_n)_{n\in P}$.
  If $P = \N$ we write $W(A)$ for $W_P(A)$, and if 
  $P = P_p = \{1, p, p^2, p^3,\ldots\}$ for a prime $p$, we write $W_p(A)$
  for $W_P(A)$.  
  We consider the Witt polynomials $w_n$ as set-theoretic maps
  $w_n: W_P(A) \to A$ for $n\in P$, and we write
  \[ w_* = (w_n)_{n\in P}: W_P(A) \To A^P. \]
  For $x\in W_P(A)$, the values $w_n(x)$ for $n\in P$ are
  called the \emph{ghost components} of $x$, and the coordinates
  $x_n$ are the \emph{Witt components}.

\end{defn}

The reason we do not write $W_P(A)$ for the codomain as
well as the domain of $w_*$ above is because they will soon have different
ring structures. 

\begin{rem} \label{rem:invertbij}
  Let $P$ be a divisor-stable set, and let $A$ be a ring such that all
  elements of $P$ have inverses in $A$.  Then we can solve for the Witt
  components of $x\in W_P(A)$ in terms of its ghost components $w_n(x)$ for
  $n\in P$, so $w_*: W_P(A) \to A^P$ is in fact a \emph{bijection}.
  Similarly, if no element of $P$ is a zero-divisor in $A$, then $w_*$ is
  an \emph{injection}.
\end{rem}

Now we can state the main theorem of this section:

\begin{thm} \label{thm:construction}
  Let $P$ be a divisor-stable set.  There is a unique covariant functor
  $W_P:\Alg_\Z\to\Alg_\Z$, such that for any ring $A$,
  \begin{enum}
  \item $W_P(A) = \prod_{n\in P} A = A^P$ as sets, and for a ring homomorphism
    $f: A\to B$, 
    \[ W_P(f)((a_n)_{n\in P}) = (f(a_n))_{n\in P}. \]
  \item The maps $w_n: W_P(A) \to A$ are homomorphisms of rings for all
    $n\in P$.
  \end{enum}
  The zero element of $W_P(A)$ is $(0,0,\ldots)$, and the unit element is
  $(1,0,0,\ldots)$. 

\end{thm}

We will give the proof of Theorem~\ref{thm:construction} in the next
section.  We will devote the rest of this section to some of its
consequences. 

\begin{defn}
  With the notation as in Theorem~\ref{thm:construction}, we call the ring
  $W_P(A)$ the 
  \emph{ring of $P$-Witt vectors, or $P$-Witt ring, with coefficients in $A$}.
  We call $W(A)$ the \emph{big Witt ring with coefficients in $A$}, and
  for a prime $p$, $W_p(A)$ is the 
  \emph{$p$-Witt ring with coefficients in $A$}.  When confusion about $P$
  is impossible, will simply call $W_P(A)$ the \emph{Witt ring} or 
  \emph{ring of Witt vectors}.

\end{defn}

\begin{rem}[Important!] \label{rem:nonstandard.notn}
  We should emphasize that the $p$-Witt ring $W_p(A)$ for a prime $p$ is
  by far the most commonly used in practice (at least in number theory).
  In fact, in the  presence of a fixed 
  prime number $p$, authors will generally write $W(A)$ for $W_p(A)$, and
  refer to $W_p(A)$ as ``the'' ring of Witt vectors over $A$.  In this
  case, people generally index the Witt components and ghost components of
  a Witt vector $x\in W_p(A)$ by the exponent of $p$: in other words,
  people write $x = (x_0,x_1,x_2,\ldots)$ for
  $(x_{p^0},x_{p^1},x_{p^2},\ldots)$ and 
  $w_n(x)$ for $w_{p^n}(x)$, $n\geq 0$.  The notation we use in this paper is
  therefore \emph{nonstandard} for number-theoretic applications (although
  it is natural in our more general context).
\end{rem}

\begin{rem} \label{rem:aftermainthm}
  \begin{enumerate}
  \item If $A$ is a $K$-algebra, it is not true in general that $W_P(A)$
    is a $K$-algebra.  For example, if $A = \F_p$ and 
    $P = \{1, p, p^2, p^3,\ldots\}$ then $W_P(\F_p)\cong\Z_p$ by 
    Theorem~\ref{thm:wittpring}, which is not an $\F_p$-algebra.
    We may still consider $W_P$ as a functor from $\Alg_K$ to $\Alg_\Z$.

  \item Let $R = \Z[\{X_n~:~n\in P\}]$.  Then for any ring $A$,
    $\Hom(R,A)$ is naturally identified with $W_P(A)$ as sets, and hence
    $W_P$ is representable.  The ring structure on $W_P(A)$ makes
    $R$ into a ring object in $\Alg_\Z$.

  \item When all elements of $P$ have an inverse in $A$, the condition
    that each $w_n$ be a ring homomorphism implies by
    Remark~\ref{rem:invertbij} that $w_*:W_P(A)\isom A^P$ is a ring
    isomorphism, where $A^P$ has the product ring structure.  Similarly,
    when the elements of $P$ are not zero-divisors in $A$, the map
    $w_*$ makes $W_P(A)$ into a subring of $A^P$; however,
    $w_*(W_P(A))\neq A^P$ in general.  

  \item Witt originally thought of the rings $W_p(A)$ as inverse limits of
    the rings $W_{P_{p(n)}}(A)$, where 
    $P_{p(n)} = \{1,p,p^2,\ldots,p^n\}$.  Since the elements of 
    $W_{P_{p(n)}}(A)$ have finitely many components, Witt thought of them
    as vectors.  This is the only sense in which rings of Witt vectors are
    related to vectors; really they are rings, and in fact ring-valued
    functors. 

  \end{enumerate}

\end{rem}

At this point it is convenient to explain how
Theorem~\ref{thm:construction} relates to Section~\ref{sec:motivation}.
Let 
\[ R = \Z[\{X_n,Y_n~:~n\in P\}], \]
and let $X = (X_n)_{n\in P},~Y = (Y_n)_{n\in P}\in W_P(R)$.  Let
$S = X + Y$ be the sum in the ring $W_P(R)$.  By definition, 
\[ w_n(S) = w_n(X) + w_n(Y) \]
for $n\in P$, so the solutions $S_n$ to the polynomial equations
\[ w_n((S_n)_{n\in P}) = w_n((X_n)_{n\in P}) + w_n((Y_n)_{n\in P}) \]
are contained in $\Z[\{X_i,Y_i:i\in P\}]$.  A similar result holds for
multiplication: namely, if $Z = X\cdot Y$, then
\[ w_n(Z) = w_n(X)\cdot w_n(Y), \]
and $Z_n\in\Z[\{X_i,Y_i:i\in P\}]$.  These polynomials in fact give the
ring laws for $W_P(A)$ for any ring $A$ (and any $P$, for that matter), as
the following corollary shows:

\begin{cor} \label{cor:ringlaws}
  Let $R$, $P$, $X$, $Y$, $Z$, and $S$ be as above.  The
  polynomials $S_n$ and $Z_n$ do not depend on the choice of $P$, and in
  addition, 
  \[ S_n,Z_n\in\Z[\{X_i,Y_i~:~i\mid n\}]. \]
  Let $A$ be an arbitrary ring, and let 
  $x,y\in W_P(A)$.  
  Let $s_n = S_n$ evaluated at the $x_i$ and $y_i$, and similarly for
  $z_n$.  Then
  \[ s = (s_n)_{n\in P} = x + y \qquad\text{and}\qquad
  z = (z_n)_{n\in P} = x\cdot y, \]
  where the addition and multiplication takes place in $W_P(A)$.
\end{cor}

\begin{pf*}
  Solving the equation $w_n(S) = w_n(X) + w_n(Y)$ 
  explicitly for $S_n$ shows that $S_n$ only depends on 
  $\{X_i,Y_i~:~i\mid n\}$.  As the equation $w_n(S) = w_n(X) + w_n(Y)$ does
  not depend on $P$, neither does $S$.  The same statements hold with $Z$
  replacing $S$.

  Define a ring homomorphism $f:R\to A$ by $f(X_i) = x_i$ and 
  $f(Y_i) = y_i$.  Since $W_P(f)$ is a ring homomorphism, we have
  \[ s = W_P(f)(S) = W_P(f)(X) + W_P(f)(Y) = x + y. \]
  By the same argument, $z = x\cdot y$.
\end{pf*}

\begin{rem} \label{rem:universal.case}
  Corollary~\ref{cor:ringlaws} essentially shows that the ring laws of
  $W_P(A)$ for an arbitrary ring $A$ can be calculated in $W_P(R)$ where
  $R = \Z[\{X_n,Y_n~:~n\in P\}]$.  Since $R$ is a ring which is
  \emph{torsion-free as a $\Z$-module}, one can often prove statements about
  $W_P(R)$ using the injection $w_*: W_P(R)\inject R^P$ of
  Remark~\ref{rem:aftermainthm}(3), and then derive facts about $W_P(A)$.
  This often-used trick is called ``reduction to the universal case'', and
  it is extremely powerful 
  --- one often cares about Witt rings over rings $A$ of characteristic $p$,
  but in order to prove theorems about these rings, one reduces to the
  characteristic-$0$ case.  
  We will make extensive use of this strategy; for
  instance, in Section~\ref{sec:exp} we will define a kind of
  characteristic-$p$ exponential map. 
\end{rem}

\begin{eg}
  Let $p$ be a prime number, and let $P = P_p = \{1,p,p^2,p^3,\ldots\}$.
  Note that 
  \[ w_{p^n} = \sum_{i=0}^n p^i X_i^{p^{n-i}}
  = X_1^{p^n} + pX_p^{p^{n-1}} + \cdots + p^{n-1}X_{p^{n-1}}^p + p^n X_{p^n},
  \]
  which agrees with our previous definition of $w_{p^n}$, except that we
  have renamed the variable $X_n$ to $X_{p^n}$
  (cf. Remark~\ref{rem:nonstandard.notn}).   
  By Corollary~\ref{cor:ringlaws}, if the $S_{p^n}$ are defined
  such that 
  \[ w_{p^n}(S_1,S_p,\ldots,S_{p^n}) 
  = w_{p^n}(X_1,X_p,\ldots,X_{p^n}) + w_{p^n}(Y_1,Y_p,\ldots,Y_{p^n}) \]
  for all $n$, then $S_{p^n}\in\Z[X_{p^i},Y_{p^i}]_{i=0,1,2,\ldots,n}$, and
  similarly for $Z_{p^n}$ and multiplication, thus verifying a claim that
  we made in Section~\ref{sec:motivation}.  In this sense, for any ring $R$,
  the ring laws of $W_p(R)$ are given by the same algebra as the ring laws
  of a strict $p$-ring.  

  To be explicit, we calculated in Section~\ref{sec:motivation} that
  \[ S_1 = X_1 + Y_1 \quad\text{and}\quad
  S_p = X_p + Y_p - \sum_{n=1}^{p-1}\frac 1p\binom pn X_1^nY_1^{p-n}. \]
  As for multiplication, we have
  \[ Z_1 = w_1(Z_1) = w_1(X_1)w_1(Y_1) = X_1Y_1 \]
  and
  \[\begin{split}
    Z_1^p + pZ_p &= w_p(Z_1,Z_p) = w_p(X_1,X_p)w_p(Y_1,Y_p) 
    = (X_1^p+pX_p)(Y_1^p+pY_p) \\
    &= (X_1Y_1)^p + pX_1^pY_p + pX_pY_1^p + p^2X_pY_p \\
    \implies Z_p &= X_1^pY_p + X_pY_1^p + pX_pY_p.
  \end{split}\]

\end{eg}

Given the above example, we can show how the
Witt rings offer a construction of the strict $p$-ring with a given
perfect residue ring of characteristic $p$.

\begin{thm} \label{thm:wittpring}
  Let $K$ be a perfect ring of characteristic $p$, and let $R$ be the
  strict $p$-ring with residue ring $K$, with Teichm\"uller reprsentatives
  $\tau:K\to R$.  Then the map $f:W_p(K) \to R$ given by
  \[ f(x_1,x_p,x_{p^2},\ldots) 
  = \sum_{n=0}^\infty \tau(x_{p^n}^{1/p^n})\,p^n \]
  is a ring isomorphism.
\end{thm}

\begin{pf*}
  It is clear that $f$ is a well-defined bijection, and that
  $f(1) = \tau(1) = 1$.  
  Let $x,y\in W_p(K)$, and let $s = x + y$.
  We will show that $f(s) = f(x) + f(y)$, and
  leave the analogous proof that $f(xy) = f(x)f(y)$ to the reader.
  By Corollary~\ref{cor:ringlaws},
  \[ s_{p^i} = S_{p^i}(x_1,y_1,x_p,y_p,\ldots,x_{p^i},y_{p^i}), \]
  and hence, since $S_{p^i}$ is a polynomial with integer coefficients,
  \[ s_{p^i}^{1/p^i} =
  S_{p^i}\big(x_1^{1/p^i},y_1^{1/p^i},x_p^{1/p^i},y_p^{1/p^i},\ldots,
  x_{p^i}^{1/p^i},y_{p^i}^{1/p^i}\big). \]
  Let $\td x_j = x_{p^j}^{1/p^j}$, and similarly for $\td y_j$ and 
  $\td s_j$.  Substituting into the above equation, we have
  \[ \td s_i = S_{p^i}\big(\td x_0^{1/p^i},\td y_0^{1/p^i},
  \td x_1^{1/p^{i-1}},\td y_1^{1/p^{i-1}},\ldots,
  \td x_{p^i},\td y_{p^i}\big), \]
  so by Theorem~\ref{thm:wittpadd}, 
  \[ \sum_{i=0}^n p^i \tau(x_{p^i}^{1/p^i}) +
  \sum_{i=0}^n p^i \tau(y_{p^i}^{1/p^i})
  \equiv \sum_{i=0}^n p^i \tau(s_{p^i}^{1/p^i}) \pmod{p^{n+1}}. \]
  Thus $f(x+y) \equiv f(x) + f(y)\pmod{p^{n+1}}$ for any $n$, completing
  the proof.

\end{pf*}

We will give another proof of Theorem~\ref{thm:wittpring} using Frobenius
and Verschiebung maps in Section~\ref{sec:fv}.

\begin{rem}
  One amazing aspect of the universal construction of $W_P(A)$ is that,
  not only can we recover the standard generalization of the ring 
  $\Z_p = W_p(\F_p)$
  to any perfect residue field $K$ of characteristic $p$ (namely, the
  unique unramified extension of $\Z_p$ with residue field $K$), but we
  can in fact define the ring $W_p(A)$ for \emph{any} ring $A$, of
  arbitrary characteristic.  In other words, the same ring laws used to
  construct $\Z_p$ can be used to define its analogue for an
  \emph{arbitrary} residue ring, although now $A = W_p(A)/\ker(w_1)$
  instead of $W_p(A)/p$.  Perhaps even more amazingly, these Witt rings
  will come equipped with Frobenius and Verschiebung maps associated to
  numbers $n\in P$ which are not necessarily prime;
  cf. Section~\ref{sec:fv}. 
\end{rem}

To end this section, we give an immediate corollary of
Corollary~\ref{cor:ringlaws}: 

\begin{cor} \label{cor:quotients}
  Let $P$ and $P'$ be divisor-stable sets, with $P'\subset P$.  The
  quotient map
  \[ (x_n)_{n\in P} \mapsto (x_n)_{n\in P'}:~ W_P(A)\To W_{P'}(A) \]
  is a ring homomorphism for any ring $A$, and hence defines a natural
  transformation $W_P\to W_{P'}$ of ring-valued functors.
\end{cor}

\section{Proof of the Existence of the Witt Rings}
\label{sec:constructionproof} 

The following argument is based on an exercise in Lang's \emph{Algebra}
\cite{lang}.  The exercise number varies by edition, but can
be located by looking under the ``Witt vectors'' entry in the index.  We
highly recommend that the reader do this exercise (with the caveat that
that Lang defines the Frobenius endomorphism incorrectly), but for
completeness we include our solution here.

Lang states that Witt told him the following
proof, but that it differs from his proof in \cite{witt}.  The proof below
uses the same ideas as in \cite{hazewinkel}.

\begin{defn}
  For a ring $A$, we let $\Lambda(A)$ be the (multiplicative) abelian group
  \[ \Lambda(A) = 1 + tA\ps t. \]
\end{defn}

\begin{lem} \label{lem:intcoeffs}
  Let $A$ be a ring.  Every element 
  $f = 1 + \sum_{n=1}^\infty x_n\,t^n\in\Lambda(A)$ can be written in the
  form $f = \prod_{n=1}^\infty (1-y_n\,t^n)$, for unique elements $y_n\in A$.
  Furthermore, there are polynomials $Y_n\in\Z[X_1,\ldots,X_n]$ and
  $X_n'\in\Z[Y_1',\ldots,Y_n']$, 
  independent of $A$, such that $y_n = Y_n(x_1,\ldots,x_n)$
  and $x_n = X_n'(y_1,\ldots,y_n)$.
\end{lem}

\begin{pf*}
  First we will prove by induction on $n$ that there are unique
  $y_1,\ldots,y_n$ such that 
  \[ f(t)/\prod_{i=1}^n(1-y_it^i) = 1 + O(t^{n+1}). \]
  The case $n = 0$ is
  clear.  Supposing the claim to be true for $n-1$, write
  $f(t)/\prod_{i=1}^{n-1}(1-y_it^i) = 1 + zt^n + O(t^{n+1})$.  We calculate
  that for $y\in A$,
  \[ (1 + y t^n)\inv(1+zt^n + O(t^{n+1}))
  = 1 + (z - y) t^n + O(t^{n+1}), \]
  which has zero $t^n$-coefficient if and only if $y = z$.  Hence there is
  a unique value $y_n = z$ such that 
  $f/\prod_{i=1}^n(1-y_it^i) = 1 + O(t^{n+1})$.

  The above proof shows that there are unique $y_1,y_2,\ldots\in A$ such
  that 
  \[ f(t) = \prod_{n=1}^\infty (1-y_n\,t^n). \]
  To show that $y_n\in\Z[X_1,\ldots,X_n]$, choose 
  $A = \Z[X_1,X_2,\ldots]$ and $f = 1 + \sum_{n\geq 1} X_n\,t^n$.  In this
  case we necessarily have $Y_n\in A$, and it is easy to see that $Y_n$
  only depends on the first $n+1$ terms of $f$, so in fact
  $Y_n\in\Z[X_1,X_2,\ldots,X_n]$.  For an arbitrary ring $B$, define a
  ring homomorphism $A\to B$ by $X_n\mapsto x_n$ for some choice of
  $x_n\in A$, and extend to a map
  $\Lambda(A)\to\Lambda(B)$.  Taking the image of the equality
  $\prod (1-Y_nt^n) = 1 + \sum X_nt^n$ in $\Lambda(B)$, we find that
  $\prod (1-y_nt^n) = 1 + \sum x_nt^n$, where 
  $y_n = Y_n(x_1,\ldots,x_n)$. 

  It is clear that the $x_n$ are integer polynomials in the $y_n$.

\end{pf*}

\begin{cor}
  For any ring $A$, the map 
  $x \mapsto f_x: W(A)\to\Lambda(A)$
  defined by
  \[ f_x(t) = \prod_{n=1}^\infty(1-x_n\,t^n)
  \quad\text{where}\quad x = (x_1,x_2,x_3,\ldots) \]
  is a bijection.
\end{cor}

Let $A$ be a $\Q$-algebra.  The Mercator series defines a
bijection $\log:\Lambda(A)\isom tA\ps t$, whose inverse is given by the
standard exponential series $\exp:tA\ps t\isom\Lambda(A)$.  Of course the
$\log$ map takes products to sums, as this is a formal property of the
Mercator series, so $\log$ is an isomorphism of abelian groups.  It is
clear that $f\mapsto -tdf/dt: tA\ps t\isom tA\ps t$ is also an isomorphism
of abelian groups, with inverse $\int -t\inv(\cdot)\,dt$.  Set
\[ D = -t\frac d{dt}\log: \Lambda(A) \isom tA\ps t. \]

\begin{lem} \label{lem:gfforw}
  Let $A$ be a $\Q$-algebra, and 
  let $x\in W(A)$.  Then
  \[ D(f_x(t)) = \sum_{n=1}^\infty w_n(x)\,t^n. \]
\end{lem}

\begin{pf*}
  The logarithmic derivative satisfies the standard identity
  \[ \frac d{dt}\log f(t) = \frac{f'(t)}{f(t)}. \]
  It is not hard to see that
  \[ \log\left(\prod_{n=1}^\infty(1-x_n\,t^n)\right) 
  = \sum_{n=1}^\infty \log(1-x_n\,t^n), \]
  so that 
  \[\begin{split}
    D(f_x(t)) = -t\frac d{dt}\log(f_x(t))
    &= \sum_{n=1}^\infty \frac{nx_nt^n}{1-x_n\,t^n} \\
    &= \sum_{n=1}^\infty (nx_nt^n + nx_n^2t^{2n} + nx_n^3t^{3n} +  \cdots). 
  \end{split}\]  
  Hence the $t^n$ term of $D(f_x(t))$ is exactly
  $\sum_{d\mid n} dx_d^{n/d} = w_n(x),$
  proving the lemma.

\end{pf*}

We need one last lemma before beginning the proof of
Theorem~\ref{thm:construction}. 

\begin{lem} \label{lem:whatisprod}
  Let $A$ be a $\Q$-algebra, and let
  $x,y\in W(A)$.  Let
  \[ f(t) = \prod_{d,e\in\N} \big(1 - x_d^{m/d} y_e^{m/e} t^m\big)^{de/m}, 
  \quad\text{where $m = \lcm(d,e)$.} \]
  Then
  \[ D(f(t)) = \sum_{n=1}^\infty w_n(x)w_n(y)\,t^n. \]
\end{lem}

\begin{pf*}
  As in the proof of Lemma~\ref{lem:gfforw}, we have
  \[\begin{split}
    D(f(t)) &= \sum_{d,e\in\N} \frac{de}m D(1 - x_d^{m/d} y_e^{m/e} t^m) \\
    &= \sum_{d,e\in\N} de\frac{x_d^{m/d} y_e^{m/e} t^m}
    {1 - x_d^{m/d} y_e^{m/e} t^m} \\
    &= \sum_{d,e\in\N} de\sum_{n\in\N} (x_d^{m/d} y_e^{m/e} t^m)^n. 
  \end{split}\]
  The $t^n$-term of the above series is 
  \[ \sum_{m\mid n} \sum_{\lcm(d,e)=m} (dx_d^{n/d})(ey_e^{n/e})
  = \left(\sum_{d\mid n} dx_d^{n/d}\right)\left(\sum_{e\mid n} ey_e^{n/e}\right)
  = w_n(x)w_n(y). \]

\end{pf*}

In the proof of Theorem~\ref{thm:construction} we will argue by reduction to
the universal case, as in Remark~\ref{rem:universal.case}.  

\begin{pf}[of Theorem~\ref{thm:construction}]
  We will show that the big Witt functor $W$ exists.
  Let $A$ be a $\Q$-algebra.  By Remark~\ref{rem:invertbij}, there is a
  unique ring structure on $W(A)$ making $w_*:W(A)\isom A^\N$
  into a ring homomorphism, where the codomain has the product ring
  structure.  Since $w_*(0) = (0,0,\ldots)$ and 
  $w_*(1,0,0,\ldots) = (1,1,\ldots)$, the zero element of $W(A)$ is 
  $(0,0,\ldots)$ and the unit element is $(1,0,0,\ldots)$.  As this
  construction is obviously functorial in $A$, we have proved
  that $W$ exists and is unique on the category $\Alg_\Q\subset\Alg_\Z$.
  We must show that the ring laws are in fact defined over the integers.
  
  Let $R = \Q[X_1,Y_1,X_2,Y_2,\ldots]$, where the
  $X_i$ and $Y_i$ are indeterminates.  Let
  \[ X = (X_1,X_2,\ldots),~
  Y = (Y_1,Y_2,\ldots)\in W(R). \]
  Let $S = (S_1,S_2,\ldots)\in W(R)$ be such that 
  $f_X(t)f_Y(t) = f_S(t)$, i.e.,
  \[ \prod_{n=1}^\infty (1-X_nt^n)
  \cdot\prod_{n=1}^\infty (1-Y_nt^n) = \prod_{n=1}^\infty (1-S_nt^n). \]
  By Lemma~\ref{lem:intcoeffs} we have
  $S_n\in\Z[X_1,Y_1,X_2,Y_2,\ldots]$, and
  by Lemma~\ref{lem:gfforw},
  \[\begin{split} \sum_{n=1}^\infty w_n(S)\,t^n &= D(f_S(t))
  = D(f_X(t)f_Y(t)) = D(f_X(t)) + D(f_Y(t)) \\
  &= \sum_{n=1}^\infty (w_n(X) + w_n(Y))\,t^n. \end{split}\]
  Hence $w_*(S) = w_*(X) + w_*(Y)$, 
  so $S = X + Y$ in $W(R)$. 
  Now let $Z = (Z_1,Z_2,\ldots)\in W(R)$ be such that
  \[ f_Z(t) = \prod_{d,e\in\N} \big(1 - X_d^{m/d} Y_e^{m/e} t^m\big)^{de/m}, 
  \quad\text{where $m = \lcm(d,e)$.} \]
  Then by Lemma~\ref{lem:intcoeffs} we have
  $Z_n\in\Z[X_1,Y_1,X_2,Y_2,\ldots]$, and by
  Lemma~\ref{lem:whatisprod}, 
  \[ \sum_{n=1}^\infty w_n(Z)\,t^n = D(f_Z(t)) 
  = \sum_{n\geq 1} w_n(X)w_n(Y)\,t^n, \]
  so $w_*(Z) = w_*(X)w_*(Y)$, and hence 
  $Z = X\cdot Y$ in $W(R)$.

  Let $A$ be an arbitrary ring, and let
  $x,y\in W(A)$.  Define
  $s = x + y$ by $s_n = S_n(x,y)$, and 
  $z = x\cdot y$ by $z_n = Z_n(x,y)$.  Reasoning as in
  Corollary~\ref{cor:ringlaws}, it is clear that when $A$ is a
  $\Q$-algebra, these recover the ring laws on $W(A)$.  In any case, we
  have constructed well-defined addition and multiplication maps on 
  $W(A)$, which are functorial in $A$.  We have not yet shown that $W(A)$
  is a ring when equipped with these addition and multiplication laws.

  Suppose that $A$ embeds into a $\Q$-algebra $A'$, which is to say, that
  $A$ is torsionfree as a $\Z$-module.  Then the inclusion 
  $W(A)\inject W(A')$ respects addition and multiplication, i.e., $W(A)$ is a
  \emph{subring} of $W(A')$.  Hence $W(A)$ is a ring for such $A$.  Now
  let $B$ be an arbitrary ring, and choose a set $\{x_i\}_{i\in I}$ of
  generators of $B$ as a $\Z$-algebra.  Set $A = \Z[\{X_i\}_{i\in I}]$,
  and let $\phi:A\surject B$ be the surjective ring homomorphism such that
  $\phi(X_i) = x_i$.  Then $W(\phi): W(A)\surject W(B)$ also respects the
  addition and multiplication laws, which is to say that $W(B)$ is a
  \emph{quotient} ring of $W(A)$.  As $A$ is a torsionfree $\Z$-module,
  $W(A)$ is a ring, so $W(B)$ is a ring.

  This completes the construction of a functor $W$ satisfying the
  properties of Theorem~\ref{thm:construction}.  The unicity of the ring
  structure on $W(A)$ is proved in the same way as the previous paragraph:
  namely, we know that $W(A)$ has only one ring structure such that $w_*$
  is a ring homomorphism when $A$ is a $\Q$-algebra; hence it is determined
  when $A$ embeds into a $\Q$-algebra, and therefore, when $A$ is the
  quotient of a ring embedding into a $\Q$-algebra.

  We leave it to the reader to construct the functor $W_P$ for an
  arbitrary divisor-stable set $P$, using the same addition and
  multiplication polynomials $S_n,Z_n$ above.

\end{pf}

\section{The Standard Topology on the Witt Rings}
\label{sec:topology}

There is a natural inverse limit topology on the Witt rings, which is an
important piece of structure, as almost all maps between Witt rings that
we will see are continuous.  This topology allows one to make sense of
infinite sums of Witt vectors, which will be very useful in the sequel.
We will assume that the reader is familiar with the theory of linear
topological rings; cf. \cite[\S0.7]{ega}.

\begin{notn}
  Let $P$ be a divisor-stable set.  For $n\in\N$ write
  \[ P(n) = \{m\in P~:~m\leq n\}. \]
  It is clear that $P(n)$ is a divisor-stable set.
  Let $\pi_n = \pi_{P,n}: W_P\to W_{P(n)}$ be the projection.
\end{notn}

It is obvious from the definitions that for any $P$, 
\[ W_P(A) = \invlimm_{n\in\N} W_{P(n)}(A) \]
as rings, under the maps $\pi_n$.

\begin{defn}
  Let $P$ be a divisor-stable set, and let $A$ be a ring, equipped with
  the discrete topology.
  The \emph{standard topology} on $W_P(A)$ is by definition the inverse
  limit topology on $\invlim W_{P(n)}(A)$, which is the same as the
  product topology on $W_P(A) = A^P$.

\end{defn}

The standard topology has the following properties, the proofs of which
are obvious.

\begin{prop}
  Let $P$ be a divisor-stable set, and let $A$ be a ring.
  \begin{enumerate}
  \item The standard topology makes $W_P(A)$ into a topological ring,
    i.e., the ring laws on $W_P(A)$ are continuous.
  \item The filtered set of ideals 
    $\{\ker(\pi_n)~:~n\in\N\}$ forms
    a neighborhood base of the identity in $W_P(A)$.
  \item $W_P(A)$ is complete and Hausdorff with respect to the standard
    topology. 
  \item The sequence $x^{(n)}\in W_P(A)$ is Cauchy if and only if, for all
    $m\in P$, $\pi_m(x^{(n)})$ is constant for $n\gg 0$; the sequence
    converges to $y\in W_P(A)$ if and
    only if, for all $m\in P$, $\pi_m(x^{(n)}) = \pi_m(y)$ for $n\gg 0$.
  \item The standard topology on $W_P(A)$ is discrete if and only if 
    $P$ is finite.
  \end{enumerate}

\end{prop}

All maps between Witt rings that we have defined so far are continuous:

\begin{prop}
  Let $P$ be a divisor-stable set and let $A$ be a ring.  The following
  maps are continuous:
  \begin{enumerate}
  \item $w_n: W_P(A)\to A$ for $n\in P$, where $A$ has the discrete
    topology.
  \item $w_*: W_P(A)\to A^P$, where $A^P$ has the product topology induced
    by the discrete topology on $A$.
  \item The projection $W_P(A)\to W_{P'}(A)$ for $P'\subset P$.
  \item The homomorphism $W_P(f):W_P(A)\to W_P(B)$ for a ring $B$ and a ring
    homomorphism $f:A\to B$.
  \end{enumerate}

\end{prop}

\begin{pf*}
  \begin{enumerate}
  \item Since $\ker(\pi_n)\subset \ker(w_n)$ we have that $\ker(w_n)$ is
    open. 
  \item The product topology is defined to be the finest topology such
    that a product of continuous maps is continuous.
  \item Let $\pi_{P,P'}:W_P(A)\to W_{P'}(A)$ be the projection.  Then
    $\pi_{P,P'}\inv(\ker(\pi_{P',n})) \supset \ker(\pi_{P,n})$.
  \item This is clear because $W_P(f) = \invlim W_{P(n)}(f)$.
  \end{enumerate}

\end{pf*}

As stated above, one of the advantages of having a topology on $W_P(A)$ is
the convergence of infinite sums.  In order to take advantage of this
property, we need to make a digression on simple arithmetic in the Witt
rings. 

\begin{prop} \label{prop:addwitt}
  Let $P$ be a divisor-stable set, let $A$ be a ring, and let
  $x,y\in W_P(A)$ be Witt vectors such that for all $n\in P$, either
  $x_n = 0$ or $y_n = 0$.  Let $s = x+y$.  Then $s_n = x_n + y_n$, i.e.,
  \[ s_n =
  \begin{cases}
    x_n &\quad\text{if } x_n\neq 0 \\
    y_n &\quad\text{if } y_n\neq 0. 
  \end{cases}\]

\end{prop}

\begin{pf*}
  As we are verifying the equality of polynomial equations, it suffices to
  prove the Proposition universally, i.e., we may assume that $A$ is a
  $\Q$-algebra.  In this case, $w_*:W_P(A)\isom A^P$ is an isomorphism, so
  we may check that $w_*(s) = w_*(x) + w_*(y)$, when
  $s_n = x_n + y_n$.  Indeed, since for every $n\in P$ either $x_n = 0$ or
  $y_n=0$, we have
  \[ w_n(s) = \sum_{d\mid n} d(x_d+y_d)^{n/d}
  = \sum_{d\mid n} dx_d^{n/d} + \sum_{d\mid n} dy_d^{n/d}
  = w_n(x) + w_n(y). \]

\end{pf*}

\begin{defn} \label{defn:teichmuller}
  Let $P$ be a divisor-stable set, let $A$ be a ring, and let $a\in A$.
  We write $[a]\in W_P(A)$ for the Witt vector whose first component is
  $a$, and whose other components are zero:
  \[ [a] = (a,0,0,0,\ldots). \]
  We call $[a]$ the Teichm\"uller representative for $a$, as the following
  Proposition justifies:

\end{defn}

\begin{prop} \label{prop:multwitt}
  Let $P$ be a divisor-stable set, let $A$ be a ring, and let
  $a\in A$.  For $x\in W_P(A)$, we have
  \[ [a]x = (a^n x_n)_{n\in P}. \]
  In particular, for $a,b\in A$, $[ab] = [a][b]$, so the map
  $a\mapsto [a]: A\to W_P(A)$ is multiplicative.

\end{prop}

\begin{pf*}
  Let $y = (a^nx_n)_{n\in P}\in W_P(A)$.
  As in the proof of Proposition~\ref{prop:addwitt}, we will check that 
  $w_n(y) = w_n([a]x)$ for $n\in P$.  Indeed,
  \[ w_n(y) = \sum_{d\mid n} d(a^dx_d)^{n/d}
  = a^n \sum_{d\mid n} dx_d^{n/d} = w_n([a])w_n(x) = w_n([a]x). \]

\end{pf*}

Let $P$ be a divisor-stable set and let $A$ be a ring.  For $n\in P$ and
$a\in A$, we provisionally define $V_n[a]$ to be the Witt vector whose
$n$th Witt component is $a$, and whose other components are zero.
Let $x\in W_P(A)$ be a Witt vector with only finitely many nonzero
components $x_n$.  Then it is clear from Proposition~\ref{prop:addwitt} that
\[ x = \sum_{n\in P} V_n [x_n]. \]
With topological considerations, the following stronger statement holds:

\begin{prop} \label{prop:infinitesum}
  Let $P$ be a divisor-stable set and let $A$ be a ring.  Let
  $x\in W_P(A)$.  Then 
  \[ x = \sum_{n\in P} V_n[x_n] 
  = \lim_{N\to\infty} \sum_{n\in P(N)} V_n[x_n]. \]
\end{prop}

The proof is immediate.  See also the treatment in \cite{hazewinkel}.

\section{The Frobenius and Verschiebung Maps}
\label{sec:fv}

Let $R$ be a finite unramified extension of $\Z_p$, with residue field
$K = \F_{p^n}$ and Teichm\"uller representatives $\tau:K\to R$.  It is a
standard fact from the theory of local fields that the 
quotient map $R\to K$ induces an isomorphism from the automorphism group
of $R$ over $\Z_p$ to 
$\Gal(K/\F_p) = \{1,\Frob,\Frob^2,\ldots,\Frob^{n-1}\}$, where
$\Frob(x) = x^p$ is the Frobenius map.  In other words, there is a
canonical lift $F_p$ of $\Frob$ to $R = W_p(K)$, which one can check is
given by
\begin{equation}
  \label{eq:liftfrob}
  F_p\left(\sum_{n=0}^\infty \tau(x_n)\,p^n\right) = 
  \sum_{n=0}^\infty \tau(x_n^p)\,p^n. 
\end{equation}

Whereas it is surprising that one can define a $p$-Witt ring
$W_p(A)$ for a ring $A$ of arbitrary characteristic, it is perhaps more
surprising that $W_p(A)$ always carries a canonical lift $F_p$ of the
Frobenius map on $W_p(A)/pW_p(A)$.  In fact, the big
Witt ring $W(A)$ has a commuting family of Frobenius maps $F_n$ for
any natural number $n$.  These maps, along with their cousins the
Verchiebung maps, are very important pieces of structure
of the Witt rings, so we devote an entire section to them.

The following will be very useful when comparing homomorphisms from Witt
rings to abelian groups.  Recall that we may consider $W_P$ as a functor
$\Alg_K\to\Alg_\Z$ for any ring $K$; this added flexibility will be
useful in the proof of Theorem~\ref{thm:frobenius}.

\begin{defn}
  For a divisor-stable set $P$, let $W_P^+:\Alg_\Z\to\Ab$ denote the
  functor that assigns to each ring $A$ the additive group underlying
  $W_P(A)$.  

\end{defn}

\begin{lem} \label{lem:comparefunctors}
  Let $K$ be a ring, and let $G:\Alg_K\to\Ab$ be a covariant abelian-group
  valued functor on $\Alg_K$.  We assume that $G$ is representable, in the
  sense that there is a $K$-algebra $R$ such that 
  $G \cong \Hom_K(R,\cdot)$ as set-valued functors.  Let $P$ be a 
  divisor-stable set, and let $u,v:W_P^+\to G$ be two natural 
  transformations.  Let $A_0 = K[x]$ and let $x_0 = [x]\in W_P(A_0)$
  (cf. Definition~\ref{defn:teichmuller}).  If 
  $u_{A_0}(x_0) = v_{A_0}(x_0)$, then $u = v$. 
\end{lem}

\begin{pf*}
  Replacing $u$ with $u-v$, we must show that if $u_{A_0}(x_0) = 0$ then
  $u=0$.  Precomposing $u$ with the natural quotient $W\to W_P$, we may
  assume that $P = \N$; hence we will consider $u$ as a natural
  transformation $\Lambda\to\Ab$, such that $u_{A_0}(1-xt) = 0$.  By
  universality, we have that $u_A(1-at)=0$ for all $A\in\Alg_K$ and all
  $a\in A$. 

  First we will show by induction on the degree of $f\in 1+tA[t]$ that
  $u_A(f) = 0$.  Let $f\in 1+tA[t]$ have degree $n$, and let 
  $g(t) = t^n f(1/t)$, so $g$ is monic of degree $n$.  If $g$ has a root
  $a\in A$, then since $g$ is monic, we can write $g(t) = (t-a)g_1(t)$,
  where $g_1$ is monic and $\deg(g_1) < n$.  Hence
  $f(t) = t^n g(1/t) = (1-at)f_1(t)$, where $f_1\in 1+tA[t]$ and
  $\deg(f_1) < n$.  Since $u_A$ is a homomorphism, by induction we have
  $u_A(f) = 0$.  

  Now suppose that $g$ does not have a root in $A$.
  Let $A' = A[X]/(g(X))$.  Then $A\inject A'$, and the residue $a$ of $X$
  in $A'$ is a root of $g\in A'[t]$.  By the above argument,
  $u_{A'}(f) = 0$, so since $G(A) = \Hom_K(R,A)$ injects into
  $G(A') = \Hom_K(R,A')$, we must have $u_A(f) = 0$.

  The above proves that for any ring $A$ and any
  $x\in W(A)$ such that $x_n = 0$ for $n\gg 0$, we have
  $u_A(x) = 0$; indeed, $f_x(t)$ is a polynomial (of finite degree).

  Now let $A = K[X_1,X_2,\ldots]$, and let 
  $X = (X_1,X_2,\ldots)\in W(A)$.  By a universality argument, it suffices
  to show that $u_A(X) = 0$.  Let $\phi = u_A(X)\in\Hom_K(R,A)$, and let
  $\psi\in\Hom_K(R,A)$ correspond to the zero element of the abelian group 
  $G(A)$.  Let $x\in R$, so $\phi(x)$ and $\psi(x)$ only involve finitely
  many $X_i$.  Let $A_n = A[X_1,X_2,\ldots,X_n]\subset A$, and suppose that
  $\phi(x),\psi(x)\in A_n$.  Consider the commutative square
  \[\xymatrix{
    {W(A)} \ar[r]^{u_A} \ar[d] & {G(A)} \ar[d] \\
    {W(A_n)} \ar[r]^{u_{A_n}} & {G(A_n)} 
  }\]
  where the vertical maps are induced by the homomorphisms
  $\pi_n: A\to A_n$ given by setting $X_i = 0$ for $i > n$.  Then 
  $W(\pi_n)(X)$ is a Witt vector with finitely many nonzero components, so 
  $u_{A_n}(W(\pi_n)(X)) = 0$.  By commutivity of the square, we have 
  \[ \pi_n\circ\psi = G(\pi_n)(\psi) = 0 
  = u_{A_n}(W(\pi_n)(X)) = G(\pi_n)(u_A(X)) = \pi_n\circ\phi, \] 
  so $\pi_n(\phi(x)) = \pi_n(\psi(x))$.  But 
  $\phi(x),\psi(x)\in A_n\subset A$, so 
  \[ \phi(x) = \pi_n(\phi(x)) = \pi_n(\psi(x)) = \psi(x). \]
  As $x$ was arbitrary, this proves that $\phi = \psi$, so $u_A(X) = 0$.

\end{pf*}

\begin{eg}
  We will primarily apply Lemma~\ref{lem:comparefunctors} when
  $G = W_P^+$, which is to say, $R = K[\{X_n:n\in P\}]$; cf. 
  Remark~\ref{rem:aftermainthm}.
\end{eg}

\begin{rem}
  The preceding lemma is stated without proof in \cite{cartier}.
  Cartier claims it is true for any functor
  $G:\Alg_K\to\Ab$ taking injections to injections, a claim that I cannot
  prove. 
\end{rem}

Now we can move on to defining the Frobenius and Verschiebung maps.  As
the Verschiebung is easier to write down, we start with it.

\begin{thm} \label{thm:verschiebung}
  Let $n\in\N$, and let $P$ be a divisor-stable set.  For any ring $A$
  define $V_n:W_P^+(A)\to W_P^+(A)$ by
  $V_n((x_m)_{m\in P}) = (y_m)_{m\in P}$, where
  \[ y_m = 
  \begin{cases}
    0 &\quad\text{if } n\nmid m \\
    x_{m/n} &\quad\text{if } n\mid m.
  \end{cases}
  \]
  Then $V_n$ is a continuous homomorphism of additive topological groups,
  and in fact defines a natural transformation $W_P^+\to W_P^+$, having
  the following properties: 
  \begin{enumerate}
  \item \[ w_m(V_n(x)) = 
    \begin{cases}
      0 &\quad\text{if } n\nmid m \\
      n\cdot w_{m/n}(x) &\quad\text{if } n\mid m
    \end{cases}
    \]

  \item If $P = \N$, so $W_P(A) = W(A) \cong \Lambda(A)$, then
    $V_n:\Lambda(A)\to\Lambda(A)$ is given by
    $V_n(f(t)) = f(t^n)$.  

  \end{enumerate}

\end{thm}

\begin{pf}
  Clearly $V_n$ is a natural transformation.
  It is obvious from the definition that when $P = \N$ we have
  $V_n(f(t)) = f(t^n)$.  Hence $V_n$ is a homomorphism of additive groups
  when $P = \N$, as $(fg)(t^n) = f(t^n)g(t^n)$.  For
  arbitrary $P$ the square
  \[\xymatrix{
    {W(A)} \ar[r]^{V_n} \ar[d] & {W(A)} \ar[d] \\
    {W_P(A)} \ar[r]^{V_n} & {W_P(A)}
  }\]
  commutes, so $V_n:W_P^+(A) \to W_P^+(A)$ is a homomorphism of additive
  groups.  From the above it is clear that $V_n:W_P(A)\to W_P(A)$ is the
  inverse limit of the maps $V_n:W_{P(m)}(A)\to W_{P(m)}(A)$ for 
  $m\in\N$, so $V_n$ is continuous.

  It remains to calculate $w_m(V_n(x))$.  Let $x\in W_P(A)$ and
  $y = V_n(x)$.  If $n\nmid m$ then
  \[ w_m(V_n(x)) = \sum_{d\mid m} dy_d^{m/d} = 0 \]
  since $n\nmid d$ when $d\mid m$.  If $m = m'n$ then
  \[ w_m(V_n(x)) = \sum_{d\mid m} dy_d^{m/d}
  = \sum_{d\mid m'} (nd) y_{nd}^{m/(nd)}
  = n\sum_{d\mid m'} d x_d ^{m'/d}
  = nw_{m'}(x). \]

\end{pf}

Note that, for $x\in A$, the above definition of $V_n[x]$ coincides with
the provisional one given in Section~\ref{sec:topology}.

\begin{eg}
  When $P = P_p = \{1,p,p^2,\ldots\}$ we have
  \[ V_p(x_1,x_p,x_{p^2},\ldots) = (0,x_1,x_p,x_{p^2},\ldots) \] 
  and 
  \[ w_*(V_p(x)) = (0, pw_1(x), pw_p(x), pw_{p^2}(x),\ldots). \]
\end{eg}

\begin{thm} \label{thm:frobenius}
  Let $n\in\N$, and let $P$ be a divisor-stable set such that
  $n\cdot P = \{np~:~p\in P\}\subset P$.  There is a
  unique natural transformation of ring-valued functors 
  $F_n:W_P\to W_P$ such that for every ring $A$ and every $x\in W_P(A)$
  and $m\in P$, we have
  \[ w_m(F_n(x)) = w_{mn}(x). \]
  The map $F_n$ is continuous.
  When $P = \N$, so $W_P(A) = W(A) \cong \Lambda(A)$, the map
  $F_n:\Lambda(A)\to\Lambda(A)$ is defined by
  \[ F_n(f)(t^n) = \big(\Norm(f)\big)(t), \]
  where $\Norm = \Norm_{A\ps t/A\ps{t^n}}: A\ps t\to A\ps{t^n}$ is the
  norm map.  

\end{thm}

I believe it was Cartier in~\cite{cartier} who first realized that the
Frobenius and the norm map are related.

\begin{pf}
  First suppose that $P = \N$.  Let $A$ be a $\Q$-algebra.  The map
  \[ (x_1,x_2,x_3,\ldots) \mapsto (x_n,x_{2n},x_{3n},\ldots):  A^\N\to A^\N \]
  is a ring homomorphism, so since $w_*:W(A)\isom A^\N$ is an isomorphism,
  it is clear that there is a unique ring homomorphism 
  $F_n:W(A)\to W(A)$ satisfying the desired property.  We claim that, on
  $\Lambda(A)$, we have $F_n(f)(t^n) = \Norm(f)(t)$.  We use
  Lemma~\ref{lem:comparefunctors}, as applied to $\Alg_\Q$, to reduce the
  claim to proving that $F_n(f)(t^n) = \Norm(f)(t)$ when $f = 1-at$
  --- note that $F_n$ and $\Norm$ are both natural transformations of
  abelian-group valued functors $\Lambda(A)\to\Lambda(A)$.
  By definition we have $w_{mn}(f) = a^{mn}$, so $w_m(F_n(f)) = a^{mn}$.
  As $w_m(1-a^nt) = a^{mn}$, we must have $F_n(f)(t) = 1-a^nt$, so
  $F_n(f)(t^n) = 1-a^nt^n$.  To calculate $\Norm(f)(t)$, we choose the
  basis $1,t,t^2,\ldots,t^{n-1}$ of $A\ps t$ over $A\ps{t^n}$.  With
  respect to this basis, the matrix for multiplication by $1-at$ is
  \[ \mu_{1-at} =
  \begin{bmatrix}
    1 & -a &  0 &  0 & \cdots & 0 & 0 \\
    0 &  1 & -a &  0 & \cdots & 0 & 0 \\
    0 &  0 &  1 & -a & \cdots & 0 & 0 \\
    \vdots & &  &    & \ddots &   &   \\
    -at^n & 0 & 0 & 0 & \cdots & 0 & 1
  \end{bmatrix}\]
  Expanding about the first column, we calculate
  \[ \Norm(1-at) = \det(\mu_{1-at}) = 1 + (-1)^nat^n \cdot (-1)^{n-1}a^{n-1}
  = 1 - a^nt^n, \]
  which proves the claim.

  At this point we can \emph{define} the natural transformation
  $F_n: W\to W$ on $\Alg_\Z$ by the formula $F_n(f)(t^n) = \Norm(f)(t)$,
  as $\Norm(f)(t)$ is defined for any ring $A$.  It remains to show that
  $F_n:W(A)\to W(A)$ is a ring homomorphism for all $A$, and that $F_n$ is
  unique.  Showing that $F_n$ is a ring homomorphism is a standard universal
  argument: one just shows that the polynomials defining $F_n$
  on the ring $\Z[X_1,X_2,\ldots]$ commute with the polynomials for
  addition and multiplication, by reducing to the case of a $\Q$-algebra.
  As usual, unicity follows by functoriality: clearly $F_n$ is
  determined for any ring $A$ embedding into a $\Q$-algebra, and hence for
  any quotient of such a ring, which is to say, any ring.
  
  We must show that $F_n:W_P\to W_P$ exists and is unique for
  general $P$ satisfying $nP\subset P$.  Directly solving the polynomial
  equations $w_{mn}(x_1,x_2,\ldots) = w_m(y_1,y_2,\ldots) = w_m(F_n(x))$
  shows that $y_m$ only depends on the $x_i$ for $i\mid nm$.  Therefore, 
  the map $F_n:W(A)\to W(A)$ descends to the quotient $W_P(A)$.  Unicity
  follows from the same universality argument.

  It remains to show that $F_n$ is continuous.  As $W_P(A)$ is a
  first-countable topological space, it suffices to show that $F_n$ is
  sequentially continuous.  Let $x^{(m)}\in W_P(A)$ converge to a point
  $x\in W_P(A)$.  This means that for any $N\in\N$, we have
  $x^{(m)}_i = x_i$ for $i\leq N$ and $m\gg 0$.  Since each Witt component
  of $F_n(y)$ only depends on finitely many components of $y\in W_P(A)$,
  it is clear that $F_n(x^{(m)})\to F_n(x)$.

\end{pf}

\begin{rem} \label{rem:afterfrob}
  The proofs of Theorems~\ref{thm:verschiebung} and~\ref{thm:frobenius}
  show that $F_n$ and $V_n$ commute with the quotient maps
  $W_P\to W_{P'}$ for appropriate $P'\subset P$.
  Lemma~\ref{lem:comparefunctors} implies that $F_n$ and $V_n$ are
  determined by the equations
  \[ F_n(1-at) = 1-a^nt \qquad\text{and}\qquad V_n(1-at) = 1-at^n \]
  on $\Lambda(A)$ (cf. \cite{cartier}).
\end{rem} 

As formal consequences of Theorems~\ref{thm:verschiebung}
and~\ref{thm:frobenius}, we have the following Propositions:

\begin{prop} \label{prop:VFnm}
  Let $n,m\in\N$, and let $P$ be a divisor-stable set.  We have
  \[ V_n\circ V_m = V_{nm} = V_m\circ V_n. \]
  When $nP\subset P$ and $mP\subset P$ we also have
  \[ F_n\circ F_m = F_{nm} = F_m\circ F_n. \]
\end{prop}

\begin{pf*}
  The statement about the Verschiebung is obvious from the definition.  
  One way to prove the claim about the Frobenius is to argue as in the
  proof of Theorem~\ref{thm:frobenius}: on $\Lambda(A)$, we have
  \[ F_n\circ F_m(1-at) = F_n(1-a^mt) = 1-a^{mn}t = F_{mn}(1-at). \]
\end{pf*}

\begin{prop} \label{prop:FVinteract}
  Let $n\in\N$ and let $P$ be a divisor-stable set with $nP\subset P$.
  Let $A$ be any ring, and let $x,y\in W_P(A)$.
  \begin{enumerate}
  \item $F_n\circ V_n(x) = n\cdot x$.
  \item $V_n(F_n(x)y) = xV_n(y)$ (i.e., $V_n$ is ``$F_n\inv$-linear'').
  \item If $m$ is prime to $n$ then $V_m\circ F_n = F_n\circ V_m$.
  \item For $m\in\N$ we have $(V_nx)^m = n^{m-1}V_n(x^m)$.
  \end{enumerate}

\end{prop}

\begin{pf*}
  As all of the above are formal identities of polynomials, a standard
  universality argument allows us to assume that $A$ is a $\Q$-algebra. 
  In this case, $w_*:W_P(A)\isom A^P$ is an isomorphism, so we need only
  check the identities on ghost components.

  \begin{enumerate}
  \item For $m\in P$ we have
    \[ w_m(F_n(V_n(x))) = w_{mn}(V_n(x)) = nw_m(x). \]
    Hence $w_*F_nV_n(x) = w_*(nx)$, so $F_nV_n(x) = nx$.

  \item For $m\in P$ we have
    \[\begin{split} w_{nm}(V_n(F_n(x)y)) &= nw_m(F_n(x)y)
      = nw_{nm}(x)w_m(y) \\
      &= w_{nm}(x)w_{nm}(V_ny) = w_{nm}(xV_ny) \end{split}\]
    and when $m\in P$ but $n\nmid m$,
    \[ w_m(V_n(F_n(x)y)) = 0 = w_m(x)w_m(V_ny) 
    = w_m(xV_ny). \]
    Thus $w_*(V_n(F_n(x)y)) = w_*(xV_ny)$.

  \item By Proposition~\ref{prop:VFnm}, we may assume that $n$ and $m$ are
    prime.  For $r\in P$ and $x\in A$ we have
    \[ w_r(V_m(F_n(x))) =
    \begin{cases}
      0 &\quad\text{if } m\nmid r \\
      mw_{r/m}(F_n(x)) = mw_{nr/m}(x) &\quad\text{if } m\mid r.
    \end{cases}\]
    On the other hand,
    \[ w_r(F_n(V_m(x))) = w_{rn}(V_m(x)) = 
    \begin{cases}
      0 &\quad\text{if } m\nmid rn \\
      mw_{rn/m}(x) &\quad\text{if } m\mid rn.
    \end{cases}\]
    Since $m$ and $n$ are distinct primes, $m\mid rn$ if and only if
    $m\mid r$, so the assertion follows.

  \item Using (2), we calculate
    \[\begin{split} (V_nx)^m &= (V_nx)^{m-1}\cdot(V_nx) =
      V_n(F_n((V_nx)^{m-1})\cdot x) \\
      &= V_n(F_n(V_nx)^{m-1}\cdot x) 
      = V_n((nx)^{m-1}\cdot x) = n^{m-1} V_n(x^m). \end{split}\]

  \end{enumerate}

\end{pf*}

\begin{rem}
  It is \emph{not} in general true that $V_n\circ F_n = n$.  However,
  it is clear from Proposition~\ref{prop:Fliftsfrob} that when $n = p$ is
  prime and $pA = 0$, then $F_p$ and $V_p$ are commuting endomorphisms of
  $W_P(A)$, so in this case we do have $V_p\circ F_p = F_p\circ V_p = p$.

\end{rem}

The first part of Proposition~\ref{prop:FVinteract} is evidence that for a
prime number $p$, $F_p$ and $V_p$ deserve to be called Frobenius and
Verschiebung maps, respectively.  The following Proposition demonstrates
that fact beyond a doubt.

\begin{prop} \label{prop:Fliftsfrob}
  Let $p$ be a prime, and let $P$ be a divisor-stable set with 
  $pP\subset P$.  Let $A$ be a ring, and let 
  $x\in W_P(A)$ and $y = F_p(x)$.  Then 
  \begin{enumerate}
  \item $y_m\equiv x_m^p\pmod{pA}$ for $m\in P$, and
  \item $F_p(x)\equiv x^p\pmod{pW_P(A)}$.
  \end{enumerate}

\end{prop}

\begin{pf*}
  We immediately reduce both assertions to the case when $P = \N$.  For
  the first claim, we replace $A$ with $A/pA$; we must show that
  when $pA = 0$, we have $F_p(x_1,x_2,\ldots) = (x_1^p,x_2^p,\ldots)$.
  Let $A = \F_p[X_1,X_2,\ldots]$, where the $X_i$ are indeterminates; by
  the usual arguments, it suffices to show that 
  $F_p(X_1,X_2,\ldots) = (X_1^p,X_2^p,\ldots)$.
  Replacing $W(A)$ with $\Lambda(A)$, we want to show that
  $F_p(\prod(1-X_nt^n)) = \prod(1-X_n^pt^n)$.  Let $f = \prod(1-X_nt^n)$.
  As $F_p(f)(t^p) = \Norm(f)(t)$, where $\Norm$ is the norm from
  $A\ps t$ to $A\ps{t^p}$, it suffices to show that 
  $\Norm(f) = \prod(1-X_n^pt^{pn}) = f^p$.
  Since $f^p\in A\ps{t^p}$, we have
  \[ f^{p^2} = N(f^p) = N(f)^p. \]
  But $A\ps t$ is a ring of characteristic $p$ with no nilpotents, so the
  $p$th power map is injective, and hence $N(f) = f^p$, as desired.  

  The second assertion is more delicate.  Let $a\in A$, and recall that
  $[a] = (a,0,0,0,\ldots)$ denotes the Teichm\"uller representative of
  $a$.  By Remark~\ref{rem:afterfrob} and Proposition~\ref{prop:multwitt},
  we see that 
  $F_p[a] = [a^p] = [a]^p$.  Now let $n > 1$.  If $p\mid n$, by
  Propositions~\ref{prop:VFnm} and~\ref{prop:FVinteract} we have
  $(V_n[a])^p = n^{p-1}V_n([a]^p)\equiv 0\pmod p$, and
  \[ F_pV_n[a] = F_pV_pV_{n/p}[a] = pV_{n/p}[a] \equiv 0
  \equiv (V_n[a])^p \pmod p. \]
  If $p\nmid n$,
  \[ F_pV_n[a] = V_nF_p[a] = V_n[a^p] = V_n([a]^p), \]
  and by Fermat's little theorem,
  \[ (V_n[a])^p = n^{p-1}V_n([a]^p) \equiv V_n([a]^p) = F_pV_n[a]\pmod p. \]
  This shows that $F_pV_n[a]\equiv (V_n[a])^p\pmod p$ for all $n\in\N$ and
  $a\in A$, i.e., $F_p(x) \equiv x^p\pmod p$ for all Witt vectors 
  $x\in W(A)$ with only one nonzero component.  

  Let $x\in W(A)$ be an arbitrary Witt vector, and let 
  $x^{(m)} = \sum_{i=1}^m V_i[x_i]$, so $x^{(m)}$ has finitely many
  nonzero components, and $x^{(m)}\to x$ as $m\to\infty$.  By the above,
  we have 
  \[ F_p(x^{(m)}) = \sum_{i=1}^m F_pV_i[x_i] \equiv \sum_{i=1}^m (V_i[x_i])^p
  \equiv \left(\sum_{i=1}^m V_i[x_i]\right)^p = (x^{(m)})^p \pmod p. \]
  By continuity of $F_p$ and of the ring laws, 
  \[ F_p(x^{(m)}) - (x^{(m)})^p \To F_p(x) - x^p \quad\text{ as }
  m\to\infty. \]
  Assume that $A$ is a torsionfree $\Z$-module; we may do this since any
  ring is a quotient of such a ring.  Since $w_*$ is an injection, $W(A)$
  is also a torsionfree $\Z$-module, so there is a unique element
  $y^{(m)}\in W(A)$ such that $py^{(m)} = F_p(x^{(m)}) - (x^{(m)})^p$.  
  Let $W_{(n)}(A) = W_{\N(n)}(A)$ be the ring of length-$n$ Witt vectors
  with coefficients in $A$, and let $\pi_n:W(A)\to W_{(n)}(A)$ be the
  projection.  Again since $w_*$ is injective, $W_{(n)}(A)$ is torsionfree
  as a $\Z$-module.  For all $n$ we have
  \[ p\pi_n(y^{(m)}) = \pi_n(F_p(x^{(m)}) - (x^{(m)})^p), \]
  which is constant for $m\gg 0$.
  Since $p:W_{(n)}(A)\to W_{(n)}(A)$ is an injection, this shows that
  $\pi_n(y^{(m)})$ is constant for $m\gg 0$, so the $y^{(m)}$ form a
  Cauchy sequence, and hence converge to a $y\in W(A)$ such that
  $py = F_p(x) - x^p$.
  
\end{pf*}

\begin{rem}
  Neither conclusion of Proposition~\ref{prop:Fliftsfrob} is 
  true in general if we replace $p$ with an integer $n$ that is not a
  prime.  For example, the first Witt component of 
  $F_6(x)$ is 
  \[ w_6(x) = x_1^6 + 3x_3^2 + 2x_2^3 + 6x_6 \not\equiv x_1^6\pmod 6. \]
  This implies that $F_6(x)\not\equiv x^6\pmod 6$, since otherwise,
  \[ w_6(x) = w_1(F_6(x)) \equiv w_1(x)^6 = x_1^6 \pmod 6, \]
  which we just showed is false.  Similarly, the first component of 
  $F_4(x)$ is not equivalent to $x^4\pmod 4$, so
  Proposition~\ref{prop:Fliftsfrob} is even false for nontrivial prime
  powers. 

\end{rem}

Note that when $K$ is a ring of characteristic $p$,
Proposition~\ref{prop:Fliftsfrob} shows that $F_p:W_P(K)\to W_P(K)$ is
given by $F_p((x_n)_{n\in P}) = (x_n^p)_{n\in P}$, so when $K$ is perfect,
$F_p$ is an automorphism.  Comparing with \eqref{eq:liftfrob}, we see that 
when $R$ is a strict $p$-ring with residue field
$K$, then $F_p:W_p(K)\isom W_p(K)$ agrees with the natural lift of
Frobenius on $R$ under the isomorphism of Theorem~\ref{thm:wittpring}.

We can also re-prove Theorem~\ref{thm:wittpring} quite easily using the
Frobenius and Verschiebung maps.  In addition, we obtain that the
isomorphism of Theorem~\ref{thm:wittpring} is a homeomorphism with respect
to the standard topology on $W_p(K)$ and the $p$-adic topology on $R$.

\begin{thm}[re-proof of Theorem~\ref{thm:wittpring}] \label{thm:wittpring2}
  Let $K$ be a perfect ring of characteristic $p$, and let $R$ be the
  strict $p$-ring with residue ring $K$, with Teichm\"uller reprsentatives
  $\tau:K\to R$.  Then the map $f:W_p(K) \to R$ given by
  \[ f(x_1,x_p,x_{p^2},\ldots) 
  = \sum_{n=0}^\infty \tau(x_{p^n}^{1/p^n})\,p^n \]
  is an isomorphism of topological rings.
\end{thm}
\newpage
\begin{pf*}
  First we will prove that $W_p(K)$ is a strict $p$-ring with residue ring
  $K$.  For $x\in W_p(K)$ we have
  \begin{equation}
    \label{eq:multbyp}
    px = F_pV_p(x) = F_p(0,x_{p^0},x_{p^1},x_{p^2},\ldots)
    = (0,x_{p^0}^p,x_{p^1}^p,x_{p^2}^p,\ldots). 
  \end{equation}
  Since $K$ is perfect, $pW_p(K) = V_p(W_p(K)) = \ker(w_1)$.  Hence
  \[ W_p(K)/pW_p(K) \cong W_p(K)/\ker(w_1) \cong K. \]
  At this point it is also easy to see that the standard
  topology and the $p$-adic topology on $W_p(K)$ coincide, so that 
  $W_p(K)$ is $p$-adically Hausdorff and complete.  It is clear by
  \eqref{eq:multbyp} that $p$ is not a zero-divisor in $W_p(K)$.
  Hence $W_p(K)$ is a strict $p$-ring with residue ring $K$, and
  Teichm\"uller representatives 
  $t:K\to W_p(K)$ given by $t(a) = [a]$, as in
  Definition~\ref{defn:teichmuller}.  

  By Theorem~\ref{thm:strictprings}, the map
  $f:W_p(K)\to R$ given by 
  \[ f\left(\sum_{n\in\N} t(x_{p^n})\,p^n\right) = \sum \tau(x_{p^n})\,p^n \]
  is an isomorphism of rings.  Since
  $t(x_{p^n})p^n = (F_pV_p)^n[x_{p^n}] = V_{p^n}([x_{p^n}^{p^n}])$,
  we have that 
  \[ (x_{p^0},x_{p^1},x_{p^2},\ldots)
  = \sum_{n\in\N} t(x_{p^n}^{1/p^n})\,p^n \]
  which completes the proof.

\end{pf*}

\section{Almost-Universal Properties}

I am not aware of any universal property that characterizes the Witt rings
$W_P(A)$ for arbitrary $A$.  However, there are some conditions under which
there exist canonical maps to and from Witt rings.
Zink \cite{zink} attributes the maps defined in Theorem~\ref{thm:froblift}
and Corollary~\ref{cor:ww2} to Cartier.  A slightly more general version
can be found in~\cite{hazewinkel}.

Recall that $\wp(P)$ denotes the set of prime numbers in $P$.

\begin{thm} \label{thm:froblift}
  Let $P$ be a divisor-stable set, and let $A$ be a ring such that no
  element of $P$ is a zero-divisor in $A$.  Suppose that $A$ is equipped
  with ring endomorphisms $\sigma_p:A\to A$ for all $p\in\wp(P)$, such that:
  \begin{enumerate}
  \item $\sigma_p(x) \equiv x^p \pmod p$ for all $x\in A$, and
  \item $\sigma_p\circ\sigma_q = \sigma_q\circ\sigma_p$ for all
    $p,q\in\wp(P)$. 
  \end{enumerate}
  Let $n\in P$, and let $n = p_1^{e_1}\cdots p_r^{e_r}$ be its prime
  factorization.  Define 
  \[ \sigma_n = \sigma_{p_1}^{e_1}\circ\cdots\circ\sigma_{p_r}^{e_r}. \]
  Then there is a unique ring homomorphism 
  $\phi:A\to W_P(A)$ such that $w_n\circ\phi = \sigma_n$ for all 
  $n\in P$. 
\end{thm}

\begin{pf*}
  Let $T\subset\Z$ be the multiplicative subset generated by $P$, and
  let $A' = T\inv A$, so $A'$ contains $A$, and
  $w_*$ defines an isomorphism $W_P(A')\isom(A')^P$.  
  If $f:A\to A$ is a ring endomorphism then  
  $f(T)=T$, so $f$ extends uniquely to
  a ring endomorphism $f':A'\to A'$.  In particular, the endomorphisms
  $\sigma_n$ extend to $\sigma_n':A'\to A'$.  Define
  \[ \phi' = w_*\inv\circ\prod_{n\in P}\sigma_n': 
  A'\To (A')^P \To W_P(A'), \]
  so $\phi'$ is the unique ring homomorphism such that 
  $w_n\circ\phi' = \sigma_n'$ for all $n\in P$.  We need only show that
  $\phi'(A)\subset W_P(A)$, as if this were true then $\phi = \phi'|_A$
  would satisfy the required properties.

  Let $x\in A$, and let $\phi'(x) = y = (y_n)_{n\in P}$.  
  We will show by induction on $n$ that $y_n\in A$.
  When $n = 1$, we have $y_1 = w_1(y) = \sigma_1(x) = x\in A$.  
  Now let $n\in P$ be arbitrary, and 
  let $n = p_1^{e_1}\cdots p_r^{e_r}$ be the prime factorization of $n$.
  Choose a prime factor $p_i$ of $n$, and set $m = n/p_i$, so by induction,
  $\sigma_m(y) = \sum_{d\mid m} dy_d^{m/d}$ with $y_d\in A$.  Therefore,
  \[ \sigma_n(y) = \sigma_{p_i}(y) \circ \sigma_m(y)
  = \sigma_{p_i}\left(\sum_{d\mid m} dy_d^{m/d}\right)
  = \sum_{d\mid m} d\sigma_{p_i}(y_d)^{m/d}. \]
  Let $d$ divide $m$, and let $p^s$ be the highest power of $p$ dividing
  $d$, so that $p^{e_i-s-1}$ is the highest power of $p$ dividing $m/d$.
  Using Lemma~\ref{lem:pcong} we have
  \[ \sigma_{p_i}(y_d)^{m/d} = 
  \big(\sigma_{p_i}(y_d)^{p_i^{e_i-s-1}}\big)^{m/(dp_i^{e_i-s-1})}
  \equiv y_d^{n/d}\pmod{p_i^{e_i-s}A}, \]
  and thus $d\sigma_{p_i}(y_d)^{m/d}\equiv dy_d^{n/d}\pmod{p^{e_i}A}$. 
  Hence 
  \[ \sigma_n(y) \equiv \sum_{d\mid m} dy_d^{n/d} 
  \equiv \sum_{d\mid n,d\neq n} dy_d^{n/d} \pmod{p_i^{e_i}A}. \]
  Since this is true for all $i$, the Chinese remainder theorem yields
  \[ \sigma_n(y) \equiv \sum_{d\mid n,d\neq n} dy_d^{n/d} \pmod{nA}, \]
  and therefore
  \[ y_n = \frac 1n\left(\sigma_n(y) - \sum_{d\mid n,d\neq n}
    dy_d^{n/d}\right) \in A. \]

\end{pf*}

The statement of Corollary~\ref{cor:ww2} is easier to comprehend with the
following bit of notation, suggested by Zink~\cite{zink}.

\begin{notn}
  Let $P$ be a divisor-stable set, and let $A$ be a ring.  For clarity, we
  write $$\hat w_n: W_P(W_P(A))\to W_P(A)$$ for the Witt polynomials,
  $\hat F_n$ for the Frobenius on $W_P(W_P(A))$, $\hat V_n$ for the
  Verschiebung, etc.  Elements of $W_P(W_P(A))$ will be denoted
  with a hat as well, and for $\hat x\in W_P(W_P(A))$, we write
  $x = (\hat x_n)_{n\in P} = (\hat x_{n,m})_{n,m\in P}$, where 
  $(\hat x_{n,m})_{m\in P}$ are the Witt components of $\hat x_n$. 

\end{notn}

\begin{cor} \label{cor:ww2}
  Let $P$ be a divisor-stable set such that $nP\subset P$ for all 
  $n\in P$.  There is a unique natural transformation of ring-valued
  functors on $\Alg_\Z$
  \[ \Delta: W_P \To W_P\circ W_P \quad\text{satisfying}\quad 
  \hat w_n\circ\Delta = F_n\quad\text{for all } n\in P. \]
  We also have the identity
  \[ W_P(w_n)\circ \Delta = F_n\qquad\text{for all } n\in P. \]
\end{cor}
\newpage
\begin{pf*}
  For any ring $A$, we have a commuting family 
  $\{F_n:W_P(A)\to W_P(A)~|~n\in P\}$ of ring endomorphisms of $W_P(A)$.
  By Proposition~\ref{prop:Fliftsfrob}, for $p\in\wp(P)$, $F_p$ lifts the
  Frobenius map on $W_P(A)/pW_P(A)$.
  Let $A$ be a ring that is torsionfree as a $\Z$-module.  Since
  $w_*:W_P(A)\to A^P$ is injective, we have in particular that no element
  of $P$ is a zero-divisor in $W_P(A)$.  Hence by
  Theorem~\ref{thm:froblift}, there is a unique map
  $\Delta_A:W_P(A)\to W_P(W_P(A))$ such that 
  $\hat w_n\circ\Delta = F_n$ for all $n\in P$.  
  Since $\hat w_*$ is an injection, it is easy to see that $\Delta_A$ is
  functorial in $A$.  Hence $\Delta$ exists and is unique on rings that
  inject into $\Q$-algebras. 

  Let $R = \Z[\{X_n:n\in P\}]$, let $X = (X_n)_{n\in P}\in W_P(R)$, and
  let $\hat Y = (\hat Y_{n,m})_{n,m\in P} = \Delta_R(X)$.  
  Then $\hat Y_{n,m}\in R$, i.e.,
  $\hat Y_{n,m}$ is an integer polynomial in the $X_n$.  For an arbitrary
  ring $A$ and $x\in W_P(A)$, define 
  $\Delta_A(x) = \hat y = (\hat y_{n,m})_{n,m\in P}$, where
  $\hat y_{n,m} = \hat Y_{n,m}(x)$.  By functoriality, this recovers the
  above definition of $\Delta_A$ when $A$ is a torsionfree
  $\Z$-module, so by the standard arguments, $\Delta_A$ is a ring
  homomorphism for arbitrary $A$.  As $\Delta_A$ is certainly functorial
  in $A$, we see that $\Delta$ exists and is unique.

  The identity $W_P(w_n)\circ\Delta = F_n$ is a relation of integer
  polynomials, so by a universality argument, it suffices to check on
  ghost components.  Let $A$ be a ring, let $x\in W_P(A)$, and let
  $\hat y = \Delta(x)$.  Since $w_m(F_n(x)) = w_{mn}(x)$, we want to show
  that for all $m\in P$, $w_m(W_P(w_n)(\hat y)) = w_{mn}(x)$.
  Indeed, $W_P(w_n)(\hat y) = (w_n(\hat y_m))_{m\in P}$, so 
  \[\begin{split} w_m(W_P(w_n)(\hat y)) &= \sum_{d\mid m} dw_n(\hat y_d)^{m/d}
  = w_n\left(\sum_{d\mid m} d\hat y_d^{m/d}\right) \\
  &= w_n(\hat w_m(\hat y)) = w_n(F_mx) = w_{nm}(x). \end{split}\]

\end{pf*}

\begin{rem}
  \begin{enumerate}
  \item 
    In proof of Corollary~\ref{cor:ww2}, we calculated the useful relation
    \[ w_m\circ W_P(w_n) = w_n\circ\hat w_m. \]

  \item Let $P$ be as in Corollary~\ref{cor:ww2}, and let $P'$ be a
    divisor-stable set containing $P$.  Then by the same proof as
    Corollary~\ref{cor:ww2}, there is a unique natural transformation
    $\Delta:W_{P'}\to W_P\circ W_{P'}$ such that
    $\hat w_n\circ\Delta = F_n$ for $n\in P$.

  \end{enumerate}

\end{rem}

There is also a bona fide universal property of the Witt rings
$W_p(K)$, where $p$ is prime and $K$ is a perfect
ring of characteristic $p$; however, I
prefer to think of it as an ``almost-universal property'' because it only
works in such limited circumstances --- really it is a property of strict
$p$-rings. 
The following definition can be found in \cite[{\S}II.5]{serre}.

\begin{defn}
  A \emph{$p$-ring} is a ring $R$ that is Hausdorff and complete for the
  topology defined by a decreasing sequence
  $\aa_1\supset\aa_2\supset\cdots$ of ideals such that
  $\aa_n\cdot\aa_m\subset\aa_{n+m}$, and such that the residue
  ring $L = R/\aa_1$ is perfect of characteristic $p$.  
\end{defn}

Theorem~\ref{thm:lameup} is a characterization of $W_p(K)$ as the
universally repelling object in the category of $p$-rings whose residue
ring is a $K$-algebra.

\begin{thm} \label{thm:lameup}
  Let $p$ be prime, and let $K$ be a perfect ring of characteristic $p$.
  Let $R$ be a $p$-ring with residue ring $L$, and let
  $f:K\to L$ be a ring homomorphism.  Then there is a unique continuous
  homomorphism $F:W_p(K)\to R$ making the square
  \begin{equation}\label{eq:upsquare}\xymatrix{
    {W_p(K)} \ar[r]^(.6){F} \ar[d]^{w_1} & R \ar[d] \\
    K \ar[r]^f & L
  }\end{equation}
  commute.

\end{thm}

\begin{pf*}
  By \cite{serre}, Chapter~II, Proposition~8, there is a unique
  multiplicative system of representatives $\tau:L\to R$.  Recall from
  Theorem~\ref{thm:wittpring2} and its proof that every element $x\in
  W_p(K)$ can be written 
  \[ x = \sum_{n=0}^\infty [x_{p^n}^{1/p^n}]\,p^n. \]
  Since power series in $p$ converge in $R$, we may define 
  $F:W_p(K)\to R$ by
  \[ F(x) = \sum_{n=0}^\infty \tau(f(x_{p^n}^{1/p^n}))\,p^n. \]
  It is clear that the square \eqref{eq:upsquare} commutes.  The proof
  that $F$ is a ring homomorphism carries over from the proof of
  Theorem~\ref{thm:wittpring} with little modification, replacing
  equivalences modulo $p^n$ with equivalences modulo $\aa_n$.
  (Alternatively, cf. \cite{serre}, Chapter~II, Proposition~9).  Since the
  standard topology on $W_p(K)$ coincides with the $p$-adic topology, and
  since $p\in\aa_1$, we see that $F$ is continuous.

  Uniqueness is proved as follows.  Let $F':W_p(K)\to R$ be another
  homomorphism satisfying the conclusions of the theorem.
  By Proposition~8 in \cite{serre}, an
  element $y\in R$ is in $\tau(L)$ if and only if $y$ is a $p^n$th power
  for all $n$.  Since $K$ is perfect, $[a]\in W_p(K)$ is a $p^n$th power
  for all $n$, so $F'([a]) = \tau(f(a))$.  Hence 
  $p^nF'([a]) = p^n\tau(f(a))$ for all $n$, so by continuity, $F = F'$.

\end{pf*}

\section{The Artin-Hasse Exponential}
\label{sec:exp}

For any prime $p$ and any ring $A$ there is a natural quotient map
$W(A)\to W_p(A)$.  It is natural to ask if that map has a section, i.e.,
if there is a natural inclusion of $W_p(A)$ into $W(A)$.  It is too much
to expect that such a section would be a ring homomorphism, but it turns
out to be true that for $\Z_{(p)}$-algebras $A$, there is a natural
homomorphism 
of abelian groups $\iota_p: W_p(A)\inject W(A)$ splitting $W(A)\surject W_p(A)$.
Zink~\cite{zink} attributes the map $\iota_p$ to Cartier.  
We will define $\iota_p$, and show how it is a kind of $p$-adic analogue of
an exponential map, which is interesting since is difficult to make sense
of the ordinary exponential series $\exp(x) = \sum_{n=0}^\infty x^n/n!$
over a ring in which there exist nonzero integers that are zero divisors.
Cartier uses $\iota_p$ in his theory of modules classifying $p$-divisible
groups, in which certain modules over the rings $W_p(A)$ are a kind of
``linearization at the origin'' of a $p$-divisible group (like a tangent
space, or a Jet space); 
in this philosophy, $\iota_p$ is in a sense the analogue of the exponential
map $\operatorname{Lie}(G)\to G$, where $G$ is a Lie group.  Cartier
theory, as well as a more thorough treatment of the Artin-Hasse
exponential, can be found in~\cite{hazewinkel}.

Cartier's construction rests on the following amazing power series:

\begin{defn}
  The Artin-Hasse exponential power series is defined by
  \[ \hexp(x) = \exp\left(x + \frac{x^p}p + \frac{x^{p^2}}{p^2} + 
    \frac{x^{p^3}}{p^3} + \cdots\right) \in \Q\ps x. \]
\end{defn}

\begin{thm} \label{thm:artinhasse}
  The Artin-Hasse exponential has $p$-integral coefficients, i.e.,
  \[ \hexp(x)\in\Z_{(p)}\ps x, \]
  where $\Z_{(p)}$ is the ring $\Z$ 
  localized at the prime ideal $(p) = p\Z$.
\end{thm}

The above theorem can be proved in many ways, including:
\begin{enumerate}
\item Dwork's criterion states that a power series 
  $f\in 1 + x\Q\ps x$ has $p$-integral coefficients if and only if 
  $f(x^p)/f(x)^p\in\Z_{(p)}\ps x$ and $f(x^p)/f(x)^p\equiv 1\pmod p$. 

\item The coefficient of $x^n$ in $n!\hexp(x)$ is the number of elements
  of the symmetric group on $n$ letters whose order is a power of $p$; the
  result then follows from a general (rather difficult) group theory fact.
\end{enumerate}
See \cite{robert} for details.  We will prove
Theorem~\ref{thm:artinhasse} in a different way, suggested by Wikipedia.

\begin{lem} \label{lem:productexpansion}
  We have the following identity in $\Q\ps x$:
  \[ \hexp(x) = \prod_{\substack{n\in\N\\p\nmid n}} (1-x^n)^{-\mu(n)/n}, \]
  where $\mu$ is the M\"obius function.
\end{lem}

It is worth clarifying that for $f\in 1+x\Q\ps x$, we define
a fractional power of $f$ as
\[ f^{a/b} = \exp(a/b\cdot\log(f)). \]

\begin{pf}
  Taking the logarithm of both sides, we want to show that 
  \[ \frac{x^p}{p} + \frac{x^{p^2}}{p^2} + \frac{x^{p^3}}{p^3} + \cdots
  = -\sum_{p\nmid n} \frac{\mu(n)}n \log(1-x^n)
  = \sum_{p\nmid n} \frac{\mu(n)}n \sum_{m=1}^\infty \frac{x^{mn}}m. \]
  The $x^r$-coefficient of the right-hand side of the above equation is
  \[ \frac 1r\sum_{\substack{n\mid r\\p\nmid n}} \mu(n). \]
  It is clear from the above equation that the $x^{p^n}$ coefficient is 
  $1/p^n$.  Let $r = p^ns$, where $p\nmid s$.  Then the $x^r$-coefficient
  is $r\inv\sum_{n\mid s} \mu(n)$, which we claim is equal to zero.  As
  $\mu(n)=0$ when $n$ is not squarefree, we may assume
  $s = p_1\cdots p_m$, where the $p_i$ are distinct primes.  Then 
  \[ \sum_{n\mid s} \mu(n) = \sum_{I\subset\{1,2,\ldots,m\}} (-1)^{\#I},
  \]
  where the sum is taken over all distinct subsets $I$ of $\{1,2,\ldots,m\}$.
  As there are $\binom{m}{i}$ such subsets of size $i$, we have
  \[ \sum_{n\mid s} \mu(n) = \sum_{i=0}^m \binom{m}{i} (-1)^i
  = (1 - 1)^m = 0. \]
\end{pf}

\begin{lem} \label{lem:1npower}
  Let $p$ be prime, let
  $f(x)\in 1 + x\Z_{(p)}\ps x$ have $p$-integral coefficients,
  and let $g(x)\in 1+x\Q\ps x$ be such that $g(x)^n = f(x)$ for an
  integer $n$ not divisible by $p$.  Then $g$ has $p$-integral
  coefficients as well, i.e., $g(x)\in 1 + x\Z_{(p)}\ps x$.
\end{lem}
\newpage
\begin{pf*}
  Write
  \[ f(x) = 1 + \sum_{m\geq 1} a_m\,x^m \quad\text{and}\quad
  g(x) = 1 + \sum_{m\geq 1} b_m\,x^m, \]
  and let $a_0 = b_0 = 1$.
  We will show inductively that $b_m\in\Z_{(p)}$.  Clearly
  $b_1 = a_1/n\in\Z_{(p)}$.  Suppose that 
  $b_1,\ldots,b_{m-1}\in\Z_{(p)}$.  Comparing the $x^m$-coefficients of 
  the equality $f(x) = g(x)^n$, we have
  \[ a_m = \sum_{\substack{m_1+\cdots+m_n=m\\0\leq m_i}} b_{m_1}\cdots b_{m_n}
  = nb_m + \sum_{\substack{m_1+\cdots+m_n=m\\0\leq m_i< m}} b_{m_1}\cdots
  b_{m_n}, \] 
  and therefore by induction,
  \[ b_m = \frac{a_m}{n} - \frac 1n
  \sum_{\substack{m_1+\cdots+m_n=m\\0\leq m_i< m}} b_{m_1}\cdots b_{m_n}  
  \in \Z_{(p)}. \]
\end{pf*}

Using Lemma~\ref{lem:productexpansion} and applying
Lemma~\ref{lem:1npower} to $f(x) = 1-x^n$, we obtain the proof
of Theorem~\ref{thm:artinhasse}.

Our goal is to find a natural section $\iota_p: W_p^+(A)\inject W^+(A)$ of
the quotient $W^+(A)\to W_p^+(A)$ for $\Z_{(p)}$-algebras $A$, and to show
that $\iota_p$ resembles an exponential map.  The
construction follows the general strategy for Witt vectors: namely, we
will make a universal construction for $\Q$-algebras $A$ (where the
logarithm and exponential do make sense), then show that the polynomials
in our construction have $p$-integral coefficients, so that it makes sense
for $\Z_{(p)}$-algebras as well.

The first thing one might try is to define $\iota_p(x)$ to be the Witt
vector $y$ such that $y_{p^r} = x_{p^r}$ and $y_n = 0$ when $n$ is not a
$p$th power.  Whereas $\iota_p$ is certainly a section of $W(A)\to W_p(A)$,
it is unfortunately not a homomorphism of additive groups.  What one
really wants is $w_n(\iota_p(x)) = 0$ for $n$ not a $p$th power, and 
$w_{p^r}(\iota_p(x)) = w_{p^r}(x)$.

Let $A$ be a $\Q$-algebra.  Consider the isomorphism 
$D:\Lambda(A)\isom tA\ps t$ of Section~\ref{sec:constructionproof}.
Define $\epsilon_p:tA\ps t\to tA\ps t$ by
\[ \epsilon_p\left(\sum_{n\geq 1} a_n\,t^n\right) = 
\sum_{r\geq 0} a_{p^r}\,t^{p^r}; \]
i.e., $\epsilon_p$ forgets the non-$p$-power coefficients of its input.
Then $\epsilon_p$ is an endomorphism of abelian groups, so there is a
unique endomorphism of $\Lambda(A)$, which we also denote by $\epsilon_p$,
such that $D\circ\epsilon_p = \epsilon_p\circ D$.  Note that 
$\epsilon_p\circ\epsilon_p = \epsilon_p$.

\begin{lem} \label{lem:epsilonp}
  Let $A$ be a $\Q$-algebra, and let 
  $f = \prod_{n\geq 1} (1-x_nt^n) \in\Lambda(A).$  Then 
  \[ \epsilon_p(f) = \prod_{r\geq 1} \hexp\big(x_{p^r}t^{p^r}\big). \]
\end{lem}

\begin{pf*}
  We calculate
  \[\begin{split}
    D\left(\prod_{r\geq 1} \hexp(x_{p^r}t^{p^r})\right)
    &= \sum_{r\geq 1} -t\frac{d}{dt}\sum_{n\geq 1}
    \frac{x_{p^r}^{p^n}t^{p^{r+n}}}{p^n}
    = -\sum_{r\geq 1}\sum_{n\geq 1} p^rx_{p^r}^{p^n}t^{p^{r+n}} \\
    &= \sum_{m\geq 1}t^{p^m} \sum_{r+n=m} p^rx_{p^r}^{p^n}
    = \sum_{m\geq 1} w_{p^m}(x)\,t^{p^m} = \epsilon_p(D(f)).
  \end{split}\]
  The Lemma follows from the definition of $\epsilon_p$.

\end{pf*}

Using the canonical identification $W(A) \cong \Lambda(A)$, we may think
of $\epsilon_p$ as an additive endomorphism of $W(A)$.
Armed with Theorem~\ref{thm:artinhasse} and Lemma~\ref{lem:epsilonp}, we
are in a position to prove:

\begin{thm} \label{thm:iota}
  Let $p$ be a prime, let $A$ be a $\Z_{(p)}$-algebra, and let
  $\pi: W(A)\to W_p(A)$ be the projection.  Define 
  $\epsilon_p:W^+(A)\to\Lambda(A)\cong W^+(A)$ by 
  \[ \epsilon_p(x_1,x_2,x_3,\ldots) = \prod_{r\geq 1}
  \hexp\big(x_{p^r}t^{p^r}\big). \]
  Then $\epsilon_p$ is an endomorphism of abelian groups, functorial in
  $A$, satisfying:
  \begin{enumerate}
  \item $\epsilon_p\circ\epsilon_p = \epsilon_p$.
  \item If $\epsilon_p(x) = y$ then $x_{p^r} = y_{p^r}$ for all $r\geq 0$.
  \item $\pi(\epsilon_p(x)) = \pi(x)$ for all $x\in W(A)$.
  \item $w_{p^r}(\epsilon_p(x)) = w_{p^r}(x)$.
  \end{enumerate}
  Furthermore, the restriction of the canonical projection $W(A)\to W_p(A)$
  induces an isomorphism $\epsilon_p W^+(A)\isom W_p^+(A)$ of abelian
  groups.  Therefore, the group homomorphism
  \[ \iota_p: W_p^+(A) \cong \epsilon_p W^+(A) \inject W^+(A) \]
  is a section of the projection $W^+(A)\to W_p^+(A)$, satisfying
  $w_{p^r}(\iota_p(x)) = w_{p^r}(x)$.
\end{thm}

\begin{pf*}
  By Theorem~\ref{thm:artinhasse}, $\epsilon_p$ is well-defined, and by
  Lemma~\ref{lem:epsilonp} and the standard universality arguments, we 
  see that $\epsilon_p$ is an additive homomorphism such that
  $\epsilon_p\circ\epsilon_p=\epsilon_p$.  Next we claim that, if 
  $y = \epsilon_p(x)$, then $y_{p^r} = x_{p^r}$ for $r\geq 0$.  It
  suffices to check the claim when $A$ is a $\Q$-algebra.  By definition,
  $w_{p^r}(y) = w_{p^r}(x)$ for all $r$, so since the $y_{p^n}$ are
  determined by the $w_{p^r}(y) = w_{p^r}(x)$, we must have
  $y_{p^n} = x_{p^n}$, as claimed.  Thus
  $\pi(\epsilon_p(x)) = \pi(x)$, so it is obvious that 
  $\epsilon_p W^+(A)$ surjects onto $W_p^+(A)$.  Injectivity is clear
  because $\epsilon_p(x)$ only depends on $\{x_{p^r}\}_{r\geq 0}$.
\end{pf*}

What Theorem~\ref{thm:iota} says is that, given 
$(x_{p^r})_{r\geq 0}\in W_p(A)$, there is a canonical choice of the
coordinates $x_n$ where $n$ is not a power of $p$, which respects the
addition law.  Explicitly, if we set $y_{p^r} = x_{p^r}$ and $y_n = 0$
when $n$ is not a power of $p$, then $\iota_p(x) = \epsilon_p(y)$.

Now we will try to indicate in what sense $\iota_p$ is a $p$-adic analogue
of an exponential map.  Let $\mc N$ be a ring without unit, such that
every element of $\mc N$ is nilpotent.  Suppose that $\mc N$ has the
structure of $\Q$-algebra.  Then we can think of the exponential map as a
group isomorphism $\exp:\mc N\isom (1+\mc N)$, where
$1+\mc N$ is the abelian group with the law $(1+a)(1+b) = 1+(a+b+ab)$.
Now suppose that $\mc N$ is a $\Z_{(p)}$-algebra.  The
statement analogous to Lemma~\ref{lem:intcoeffs} says that every element
of $\hat\Lambda(\mc N) := 1 + \mc N[t]$ can be written uniquely as a
finite product $\prod (1-a_nt^n)$.  This gives an identification of the
additive group $\hat W^+(\mc N)$ of finite-length Witt vectors with
entries in $\mc N$, with the multiplicative group $\hat\Lambda(\mc N)$,
analogous to the identification $W^+(A)\cong\Lambda(A)$ for a ring
$A$. The finite version of Theorem~\ref{thm:iota} says that 
$\epsilon_p\hat W^+(\mc N)\isom\hat W_p^+(\mc N)$.  Consider the
composition 
\[ E_p: \hat W_p^+(\mc N)\cong\epsilon_p\hat W^+(\mc N)\inject\hat W^+(\mc N) 
= 1 + t\mc N[t] \overset{t\mapsto 1}\To (1+\mc N), \]
given explicitly by
\[ E_p(x_{p^0},x_{p^1},\ldots,x_{p^m},0,0,\ldots) 
= \prod_{n=0}^m \hexp(x_{p^n}). \]
It is clear that $E_p$ is a homomorphism of abelian groups, and one can
show that its kernel is equal to $(\Id-V_p)\hat W_p^+(\mc N)$.  This
``exponential'' map is now an isomorphism 
\[ E_p: \hat W_p^+(\mc N)/(\Id-V_p)\isom(1+\mc N), \]
defined for any nilpotent $\Z_{(p)}$-algebra $\mc N$.

The preceding discussion is a special case of an isomorphism from Zink's
theory of displays~\cite{zink}.

\bibliographystyle{amsalpha}
\bibliography{witt}

\bigskip\bigskip
\hrule
\bigskip\bigskip
\noindent
Copyright 2007, Joseph Rabinoff.

\end{document}